\PassOptionsToPackage{unicode}{hyperref}
\PassOptionsToPackage{naturalnames}{hyperref}

\documentclass[a4paper,11pt]{article}

\usepackage{amsfonts,amssymb,amsmath,amsthm}	
\usepackage[english]{babel}                 	
\usepackage[T1]{fontenc}                    	
\usepackage{setspace}      						
\usepackage{graphicx}      						
\usepackage{tikz}          						
\usepackage{enumitem}
\usepackage[permil]{overpic}       				
\usepackage{fancybox}							
\usepackage{ifthen}								
\usepackage{subcaption}
\usepackage[ plainpages = false, pdfpagelabels,
	     bookmarks,
	     bookmarksopen = true,
	     bookmarksnumbered = true,
	     breaklinks = true,
	     linktocpage,
	     pagebackref,         
	     colorlinks = true,  
         hypertexnames=false,   
	     linkcolor = green!50!black,
	     urlcolor  = blue!50!black,
	     citecolor = blue!50!black,
	     anchorcolor = blue,
	     hyperindex = true,
	     hyperfigures
	     ]{hyperref}
\usepackage{tikz}
\usepackage{tikz-3dplot}
\usepackage[numbers]{natbib}
\usepackage{fancyhdr, graphicx}
\usepackage{pgfgantt}
\usepackage{bm}

\usepackage[margin=1in]{geometry}
\usepackage{bbm}
\usepackage{lineno}
\usepackage[colorinlistoftodos,prependcaption]{todonotes}
\usepackage[capitalise]{cleveref}
\usepackage{autonum} 
\usepackage{amssymb}
\usepackage{titlesec}
\usepackage{booktabs}

\usepackage{pgfplots}
\pgfplotsset{compat=newest}
\usetikzlibrary{plotmarks}
\usetikzlibrary{arrows.meta}
\usepgfplotslibrary{patchplots}
\usepackage{grffile}
\usepackage{amsmath}
\pgfplotsset{plot coordinates/math parser=false}
\newlength\figureheight
\newlength\figurewidth
\pgfplotsset{every axis/.append style={
                    label style={font=\small},
                    tick label style={font=\footnotesize},
                    legend style={font=\tiny},
                    title style={font=\small}
                    }}
\usetikzlibrary{plotmarks}
\pgfplotsset{minor grid style={dotted,gray}} 
\pgfplotsset{every axis/.append style={thick, tick style=semithick}}
\usepgfplotslibrary{groupplots}
\pgfplotsset{xticklabel style={/pgf/number format/fixed,
        /pgf/number format/precision=5},yticklabel style={/pgf/number format/fixed,
        /pgf/number format/precision=5}}
\newcommand{\inputtikz}[1]{%
  \includegraphics{#1}
}

\newcommand{\vect}[1]{\textrm{\boldmath${#1}$}} 
\DeclareMathOperator{\diag}{diag}               
\newcommand{\matr}[1]{\mathbf{#1}}              
\newcommand{\RhsV}{\mathbf{R}_\vect{v}} 
\newcommand{\advectQ}{\mathbf{F}_\vect{q}} 
\newcommand{\advectU}{\mathbf{F}_\vect{u}} 
\newcommand{\advectV}{\mathbf{F}_\vect{v}} 
\newcommand{\Q}{\mathbf{Q}} 
\newcommand{\V}{\mathbf{V}} 
\newcommand{\SourceMatrix}{\matr{G}} 
\newcommand{\SourceVector}{\vect{g}_{\nu}} 
\newcommand{\sourceScalar}{g_{\nu}} 
\newcommand{\Vhat}{\hat{\mathbf{V}}} 
\newcommand{\norm}[1]{\Vert #1 \Vert}
\providecommand{\U}{}
\renewcommand{\U}{\mathbf{U}} 
\newcommand{\Nt}{N_\mathrm{t}}
\newcommand{\email}[1]{\href{mailto:#1}{\texttt{#1}}}

{\theoremstyle{remark} \newtheorem*{remark}{Remark}}

\newcommand{\refereeOne}[1]{#1}

\newcommand{\refereeThree}[1]{#1}

\binoppenalty=\maxdimen 
\relpenalty=\maxdimen   

\title{Macro-micro decomposition for consistent and conservative model order reduction of hyperbolic shallow water moment equations:\\
A study using POD-Galerkin and dynamical low rank approximation}
\author{
Julian Koellermeier\footnote{Bernoulli Institute, University of Groningen \email{j.koellermeier@rug.nl}} ,
Philipp Krah\footnote{Institut de Mathématiques de Marseille (I2M), Aix-Marseille Université, \email{philipp.krah@univ-amu.fr}} ,
Jonas Kusch\footnote{Institut für Mathematik, Universität Innsbruck, \email{jonas.kusch1@gmail.com}}
}

\begin{document}

\maketitle

\begin{abstract}
Geophysical flow simulations using hyperbolic shallow water moment equations require an efficient discretization of a potentially large system of PDEs, the so-called moment system.
This calls for tailored model order reduction techniques that allow for efficient and accurate simulations while guaranteeing physical properties like mass conservation.

In this paper, we develop the first model reduction for the hyperbolic shallow water moment equations and achieve mass conservation. This is accomplished using a macro-micro decomposition of the model into a macroscopic (conservative) part and a microscopic (non-conservative) part with subsequent model reduction using either POD-Galerkin or dynamical low-rank approximation only on the microscopic (non-conservative) part.
Numerical experiments showcase the performance of the new model reduction methods including high accuracy and fast computation times together with guaranteed conservation and consistency properties.
\end{abstract}


\section{Introduction}

Accurate simulation of free-surface flows is necessary for prediction of natural hazards like floods, landslides, tsunami waves, and weather forecasting \cite{christen2010ramms,courtier1988global}. However, the solution of the full incompressible Navier-Stokes equations is often too costly and simplified models like the \emph{shallow water equations} (SWE) often yield inaccurate results, since the SWE model assumes a constant velocity over the vertical axis. The recently derived \emph{hyperbolic shallow water moment equations} (HSWME) overcome this problem by allowing for polynomial velocity profiles \cite{Koellermeier2020}. The model is based on a hyperbolic regularization of \cite{kowalski2018moment} using techniques from kinetic theory \cite{Fan2016,Koellermeier2014}. The increased accuracy of the model was shown in numerical simulations, including sediment transport \cite{Garres2021}. However, the HSWME also lead to a higher computational cost and memory footprint due to additional nonlinear equations for the expansion coefficients of the polynomial velocity profile. The goal of this work is to perform an additional model reduction to improve efficiency of the model while at the same time preserving important properties like conservation of mass.

To reduce computational costs during the simulation of the HSWME we propose to project the system on a lower dimensional manifold using Galerkin methods. In this context we apply two model reduction techniques. First, the classical offline-online procedure (see for a review \cite{KunischVolkwein2002}), where important modes are calculated using a \emph{proper orthogonal decomposition} (POD) in a so-called offline phase. This allows to conduct further online computations with the generated basis at significantly reduced computational costs. Second, the online adaptive basis method, \emph{dynamical low-rank approximation} (DLRA) \cite{koch2007dynamical}, that along with other methods like the \emph{adaptive basis and adaptive sampling discrete empirical interpolation method} (AADEIM) \cite{Peherstorfer2020}, adapts the dominant modes by online updates to the solution locally in time. DLRA therefore does not require an offline computation since training is shifted to the online phase.

Most commonly in \emph{model order reduction} (MOR) a reduced basis is computed with the help of the POD, that was first introduced in \cite{Lumley1967} in the context of fluid dynamics. In combination with Galerkin projection methods, POD-Galerkin has been successfully applied in fluid dynamics \cite{LassilaManzoniQarteroniRozza2014}, including the SWE \cite{LozovskiyFarthingKeesGildin2016,LozovskiyFarthingKees2017,StefanescuNavon2013,StefanescuSanduNavon2014,StefanescuSanduNavon2015}.
Since POD-Galerkin does not preserve conservation laws, one possibility is to introduce Hamilton formulations that ensure structure-preservation when constructing \emph{reduced order models} (ROM) of SWE with the help of POD \cite{KarasozneYildizUzunca2021,KarasozenYildizUzunca2022}.
In this paper we consider the HSWME that introduce $N$ additional equations to the 1D-SWE to account for non-constant velocity profiles along the vertical axis. With the aim to reduce the HSWME to a system yielding a similar complexity as the SWE we then apply POD-Galerkin as the first approach. In contrast to existing SWE-ROM approaches \cite{LozovskiyFarthingKeesGildin2016,LozovskiyFarthingKees2017,StefanescuNavon2013,StefanescuSanduNavon2014,StefanescuSanduNavon2015,KarasozneYildizUzunca2021,KarasozenYildizUzunca2022},
that try to reduce the high dimensional state-space by Galerkin projections onto the spatial modes,
we utilize POD modes that reduce the system only in the vertical direction. Here, the idea is to replace the $N+2$ dimensional HSWME system by a $r+2$ (where $r\ll N$) dimensional system. This is achieved through a linear combination of global ansatz functions that are determined by POD. Similar ideas have been first introduced by \cite{VogeliusBabuska1981P1,VogeliusBabuska1981P2,VogeliusBabuska1981P3} and later in fluid dynamics to study the pipe flow where a coarse model featuring the dominant dynamics in the flow direction is enriched by a fine model that accounts for the additional dynamics in the transverse direction with the help of a modal expansion \cite{PerottoErnVeneziani2010,PerottoVeneziani2014,PerottoRusconiVeneziani2017,MansillaBlancoBulantDariFeijo2017}.
This so called hierarchical model order reduction was further expanded for non-linear PDEs in \cite{SmentanaOhlberger2017}.

Dynamical low-rank approximation \cite{koch2007dynamical} for matrix differential equations approximates the solution by a low-rank matrix decomposition and derives evolution equations to update the factors of the solution in time. These evolution equations are determined by minimizing the defect while restricting the evolution of the solution to the manifold of low-rank matrices. Stable time integrators for the resulting DLRA system, which are robust irrespective of the curvature of the low-rank manifold \cite{kieri2016discretized}, are the matrix projector--splitting integrator \cite{lubich2014projector} as well as the ``unconventional'' \emph{basis update \& Galerkin step} (BUG) integrator of \cite{ceruti2020unconventional}. Here, we use the BUG integrator which only evolves the solution forward in time, thereby facilitating the construction of stable spatial discretizations \cite{kusch2021stability}. Moreover, the BUG integrator enables a straightforward basis augmentation step \cite{ceruti2022rank} which simplifies the construction of rank adaptive methods \cite{ceruti2022rank,kusch2021robust,hauck2022predictor} and allows for conservation properties \cite{ceruti2022rank,einkemmer2022robust}.

Comparing computational results for POD-Galerkin and DLRA, we note that the generation of basis functions during the online computation makes dynamical low-rank approximation significantly more costly than the online computation of POD-Galerkin.
However, since DLRA does not require an offline computation in which the full model needs to be evaluated repeatedly, it offers three main advantages over POD-Galerkin: First, when the full order model is expensive to compute or cannot be stored in memory at a desirable accuracy, the offline phase of POD-Galerkin can become unfeasible to compute. Second, DLRA allows a straightforward application in settings that require a retraining of the basis with POD-Galerkin. Third, since DLRA allows to adapt the basis during the online phase it is independent of the choice of snapshots and it is beneficial in settings where a precomputed linear basis might not describe the dynamics sufficiently. In particular DLRA is superior to POD-Galerkin for transport dominated systems, where the Kolmogorov $n$-width decays slowly \cite{OhlbergerRave2015,GreifUrban2019}.

While model order reduction techniques promise a reduction of cost and memory, they can struggle to preserve important quantities of the full system such as boundary conditions or mass conservation. To guarantee preserving these properties, we propose to decompose the HSWME solution into two parts: (1) macroscopic variables: water height and momentum and (2) microscopic variables: higher-order moments of the velocity profile. Separate evolution equations for the macro and micro solution parts are derived and the model order reduction strategy is only applied to the computational and memory intensive microscopic part (2). Thereby, we facilitate imposing boundary conditions for the macroscopic water height and momentum while guaranteeing mass conservation. In addition, we achieve consistency with the underlying SWE model in case of vanishing microscopic structures.

The construction of efficient numerical method for dynamical low-rank approximation requires a carefully chosen formulation of the low-rank evolution equations and the construction of adequate numerical discretizations. Therefore, besides the introduction of a macro-micro decomposition to preserve solution invariants in model order reduction, main novelties of this work are:
\begin{itemize}
\item \emph{The derivation of efficient low-rank evolution equations for HSWME.} We derive evolution equations of low-rank factors of the HSWME which do not require the computation of full-rank matrices. Thereby, the memory footprint and computational costs to evolve the dynamical low-rank approximation are significantly reduced.
\item \emph{The construction of an efficient numerical scheme.} A classical implicit time discretization of friction terms yields prohibitive costs when evolving the low-rank approximation. We propose a splitting step to eliminate computationally costly terms in the numerical scheme, allowing for a significant reduction in runtime.
\end{itemize}

The remaining parts of this paper are structured as follows: We present the main concept of the macro-micro decomposition as well as the correction terms for model order reduction in \cref{sec:Concept}. The mathematical model HSWME is briefly explained in \cref{sec:model}. \Cref{sec:discr} covers the numerical discretizations including the implementation of the macro-micro decomposition for the new model. The conservative POD-Galerkin method is adopted for the considered model in \cref{sec:POD-Galerkin}, followed by the dynamical low-rank approximation in \cref{sec:DLRA}. Numerical results showing the performance of both methods are presented in \cref{sec:numerics} and the paper ends with a brief conclusion.

\section{Concept}
\label{sec:Concept}
\begin{figure}[htp!]
    \centering
    \includegraphics[width=0.6\linewidth]{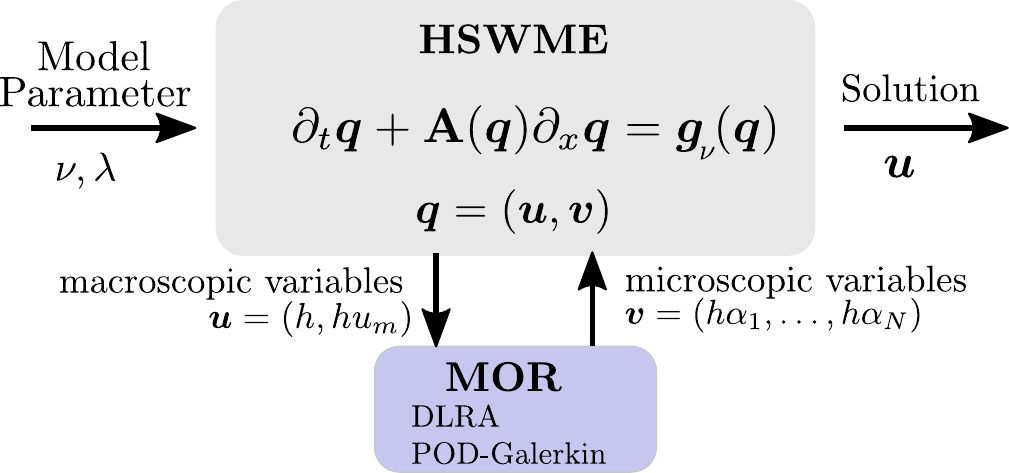}
    \caption{Concept of the macro-micro decomposition MOR approach for HSWME. The HSWME is decomposed in macroscopic variables $\vect{u}$ and microscopic variables $\vect{v}$. While the macroscopic variables are computed with a standard conservative scheme, the microscopic variables are treated by MOR techniques DLRA or POD-Galerkin. }
    \label{fig:split-cons-MOR}
\end{figure}

A core issue of model order reduction is the potential violation of crucial properties of the original problem. To preserve these properties we choose to first decompose the original dynamics into evolution equations describing the macroscopic (conservative) quantities on the one hand and microscopic (non-conservative) correction terms on the other hand. A reduced model is then derived only for the microscopic corrections terms, whereas the original macroscopic evolution equations remain unaltered. Thereby, basis functions which are important for overall solution properties remain unchanged, whereas the overall solution complexity that mainly arises from correction terms is reduced through model order reduction. The concept is visualized in \cref{fig:split-cons-MOR}.

We consider equations of the following type
\begin{equation}\label{e:eqntype}
    \partial_t \vect{q} + \matr{A}(\vect{q}) \partial_x \vect{q} = \SourceVector(\vect{q}),
\end{equation}
where $\matr{A}(\vect{q}) \partial_x \vect{q}$ denotes the transport term and $\SourceVector(\vect{q})$ denotes the friction term, which depends on a friction parameter $\nu$.
Furthermore, the variable vector $\vect{q}= (\vect{u}, \vect{v})$ contains both macroscopic variables $\vect{u}$ and microscopic variables $\vect{v}$.

For the numerical solution, we perform two main steps: (I) A first order operator splitting to \cref{e:eqntype}, in which we split the transport part from the right-hand side friction part. (II) A decomposition of the solution into macroscopic ($\vect{u}$) and microscopic ($\vect{v}$) variables. Note that the operator splitting only introduce first order errors in time and space, which are of the same order as the discretization errors in the sub-steps. The macro-micro variable decomposition does not introduce an error but a decoupled solution of the macro and micro step will do so. The procedure can be summarized by the following steps
\begin{alignat}{5}
      &\text{Step 1: transport} \qquad& \partial_t \vect{q}&+ \matr{A}(\vect{q}) \partial_x \vect{q} &\;=\;& 0\,, \label{e:transport} \\
      &\text{Step 2: friction} \qquad& \partial_t \vect{q}& &\;=\;& \SourceVector (\vect{q})\,,\label{e:friction}
\end{alignat}
and the subsequent macro-micro decomposition
\begin{alignat}{5}
      &\qquad \text{Step 1a: macro transport} & \quad\vect{q}^n= &(\vect{u}^n, \vect{v}^n) & \stackrel{\eqref{e:transport}}{\Rightarrow} (\widetilde{\vect{u}}^{n+1}, \vect{v}^n),\quad& \label{e:transport_macro}\\
      &\qquad \text{Step 1b: micro transport} & &(\widetilde{\vect{u}}^{n+1}, \vect{v}^n) & \stackrel{\eqref{e:transport}}{\Rightarrow} (\widetilde{\vect{u}}^{n+1}, \widetilde{\vect{v}}^{n+1}),& \label{e:transport_micro}\\
      &\qquad \text{Step 2a: macro friction} & \quad \widetilde{\vect{q}}^{n+1}= &(\widetilde{\vect{u}}^{n+1}, \widetilde{\vect{v}}^{n+1}) & \stackrel{\eqref{e:friction}}{\Rightarrow} (\vect{u}^{n+1}, \widetilde{\vect{v}}^{n+1}),& \label{e:friction_macro}\\
      &\qquad \text{Step 2b: micro friction} & &(\vect{u}^{n+1}, \widetilde{\vect{v}}^{n+1}) & \stackrel{\eqref{e:friction}}{\Rightarrow} (\vect{u}^{n+1}, \vect{v}^{n+1}).&  \label{e:friction_micro}
\end{alignat}

\begin{remark}
Note that macro-micro decomposition has been used in a wide variety of applications in combination with model order reduction, see e.g. \cite{einkemmer2021asymptotic,einkemmer2022asymptotic}. In contrast to these works, we employ the macro-micro decomposition not to construct asymptotic--preserving schemes, but for two purposes: 1) Guarantee mass conservation and 2) ensure consistency with the SWE when the microscopic part tends to zero. Note that the latter also ensures that the model order reduction method does not reduce the approximation accuracy of water height and mean velocity. That is, these two quantities of interest can be represented accurately, independent of any chosen basis functions. 
The choice of macroscopic and microscopic variables for a general model is not predetermined and this could change based on the interpretation of the model or a spectral analysis potentially identifying a spectral gap, e.g., \cite{Amrita2022}.
\end{remark}

\begin{remark}
    The notion of macro and micro variables originates from the limit of the SWME with respect to 
    vanishing slip length. This models a perfect slip bottom such that only the interior friction remains and drives the ``microscopic'' higher-order moments to zero resulting in a constant velocity profile \cite{Huang2022}. The remaining water height $h$ and velocity $u_m$ evolve according to the ``macroscopic'' shallow water equations. This justifies the decomposition into a macro part evolving on the equilibrium manifold and a remaining micro part. For more details, we refer to \cite{Huang2022}. 
\end{remark}

In the next sections, we will first introduce the full (coupled) model \cref{e:eqntype}, then discuss the space time discretization of the operator splitting of \cref{e:transport} and \cref{e:friction}, before deriving a detailed scheme for the macro-micro decomposition system \cref{e:transport_macro} - \cref{e:friction_micro}, which can then be efficiently used for model reduction algorithms like POD-Galerkin and DLRA of the micro steps \cref{e:transport_micro} and \cref{e:friction_micro}.

\section{Models for Shallow Flows}
\label{sec:model}
In the following, we introduce two full models for shallow water flows: the shallow water moment equations and their hyperbolic formulation, the hyperbolic shallow water moment equations.

We note that an alternative approach to starting with these models is to directly apply the model reduction techniques to the underlying incompressible Navier-Stokes equations. Conceptually, this is possible. 
In the case of DLRA, the evolution equation for the $L$-step is then continuous in the velocity variable, and a moment method can be used to obtain a discretized set of equations. This is equivalent (when using the BUG integrator as we do in this work) to first performing the moment approximation and then applying DLRA. However, using the moment method first and the model reduction second is much more convenient in terms of notation and efficiency of presentation. It also allows for a very intuitive understanding of how to carry over the desirable properties of the moment model to the numerical solution.

\subsection{Shallow water moment equations}
\label{sec:SWME}
Shallow water flows are often modeled using the standard SWE, which are hyperbolic balance laws for the macroscopic variables water height $h$ and the mean velocity $u_m$, derived from the underlying incompressible Navier-Stokes equations. For simplicity, we consider a one-dimensional vertical coordinate $x \in \mathbb{R}$, a flat bottom topography, and a Newtonian fluid. The SWE model for $h$ and $u_m$ then reads \cite{kowalski2018moment}
\begin{equation}\label{e:SWE}
    \partial_t
    \begin{pmatrix}
    h\\
    h u_m\\
    \end{pmatrix} +\partial_x
    \begin{pmatrix}
    h u_m\\
    h u_m^2 + \frac{1}{2}g h^2 \\
    \end{pmatrix} =
    -\frac{\nu}{\lambda}
   \begin{pmatrix}
    0\\
    u_m\\
    \end{pmatrix},
\end{equation}
with gravity constant $g$, slip length $\lambda$ and kinematic viscosity $\nu$. Note that the system is written in so-called \emph{convective variables} $(h, h u_m)$, which will also be the case for the more advanced shallow water Moment model below. 
The addition of a non-flat bottom topography term in \cref{e:SWE} will pose no conceptual difficulties and is left for future work, see \cite{Garres2021,Pimentel2022} for examples of the straightforward treatment.
Further note that the concept can be readily extended to the 2D case. However, so far there is very little research on 2D models with vertically resolved velocity profiles. Even though the conceptual framework of the model is laid out in \cite{kowalski2018moment} and the hyperbolic regularisation below can be extended based on \cite{Koellermeier2020}, to the best of our knowledge, no full 2D simulations exist. We therefore focus on the 1D case in the current work.

Due to its simplicity with only the macroscopic mean velocity $u_m$, the SWE model cannot represent vertical variations of the velocity, representing micro structure of the profile. The assumption of a constant velocity profile breaks down for many applications, especially when considering strong bottom friction \cite{kern_bartelt_sovilla_2010}. Even when starting with constant velocity profiles, bottom friction leads to a deceleration of the fluid close to the bottom giving rise to more complex velocity profiles.

The recently developed \emph{shallow water moment equations} (SWME) \cite{kowalski2018moment} tackle this problem by introducing a polynomial expansion of the velocity profile $u(t,x,z)$ depending on the vertical variable $z \in [0,h]$ as follows
\begin{equation}\label{e:expansion}
    u(t,x,z)=u_m(t,x)+\sum_{j=1}^{N}\alpha_j(t,x)\phi_j\left(\frac{z}{h(t,x)}\right),
\end{equation}
where $\phi_j:[0,1]\rightarrow\mathbb{R}$ are the \emph{scaled Legendre polynomials} of degree $j$ defined by
\begin{equation} \label{e:defphi}
    \phi_j(\zeta) = \frac{1}{j!} \frac{d^j}{d\zeta^j} (\zeta - \zeta^2)^j.
\end{equation}
The basis functions $\phi_j$ form an orthogonal basis, due to $\int_0^1 \phi_m(\zeta) \phi_n(\zeta) d\zeta = \frac{1}{2n+1} \delta_{m n}$,
where $\delta_{m n}$ is the Kronecker delta.

The corresponding expansion coefficients $\alpha_j$ for $j = 1,2,\ldots,N$ of the polynomial expansion \cref{e:expansion} are also called \emph{moments} in analogy to moment models from kinetic theory \cite{Torrilhon2016}. The moments represent micro variables that augment the macroscopic variables $h$ and $u_m$.

The expansion \cref{e:expansion} then allows for the representation of complex velocity profiles, e.g., linear, quadratic, ..., up to the maximal polynomial degree $N\in \mathbb{N}$. While a larger maximal degree $N$ leads to potentially higher accuracy of the velocity profile, more moments are used in the representation which leads to higher computational costs.

Evolution equations of the moments in time and space can be derived by a projection of the underlying incompressible Navier-Stokes equations onto the basis functions \cref{e:defphi}. In combination with the conservation of mass, the resulting system of equations is called the shallow water moment equations. More details on the derivation of the SWME can be found in \cite{kowalski2018moment}.

The SWME model for water height $h$, mean velocity $u_m$, and coefficients $\alpha_i$ for $i=1,\ldots,N$ can be written as
\begin{equation}\label{e:SWME_arbitrary}
    \partial_t \vect{q} + \widetilde{\matr{A}}(\vect{q}) \partial_x \vect{q} = \SourceVector(\vect{q}),
\end{equation}
with convective variable vector $\vect{q} = \left[h, h u_m, h\alpha_1, \ldots, h\alpha_N \right]^{T} \in \mathbb{R}^{N+2}$, transport matrix $\widetilde{\matr{A}}(\vect{q})\in \mathbb{R}^{{N+2} \times {N+2}}$, and right-hand side friction term $\SourceVector(\vect{q})\in \mathbb{R}^{N+2}$.

In this paper, we do not use the SWME model in the form derived in \cite{kowalski2018moment} due to its lack of hyperbolicity, which was attributed to stability issues in \cite{Koellermeier2020}. We therefore do not show the explicit form of the matrix $\widetilde{\matr{A}}(\vect{q})$ and source term $\SourceVector(\vect{q})$ of the SWME here, but directly introduce the hyperbolic regularization derived in \cite{Koellermeier2020}.

\subsection{Hyperbolic shallow water moment equations}
\label{sec:HSWME}
Based on results from kinetic theory \cite{Fan2016,Koellermeier2014,Cai2013}, the HSWME were first derived in \cite{Koellermeier2020} and have been extended and analyzed in \cite{Pimentel2022,Huang2022}. As a key property of transport dominated problems, hyperbolicity is the property of the transport matrix to have real eigenvalues corresponding to waves with real propagation speeds.  It was shown that the HSWME yield accurate results while preserving hyperbolicity of the model equations \cite{Koellermeier2020}. Recently, the model has been applied to sediment transport \cite{Garres2021} and analysis of steady states and equilibrium stability has been performed \cite{Huang2022,Pimentel2022}. 

In this paper, we use the HSWME model in its standard form written in terms of the convective variable vector $\vect{q} = [h,h u_m, h\alpha_1,\dots, h\alpha_N]^{T} \in \mathbb{R}^{N+2}$ as
\begin{equation} \label{e:HSWME}
    \partial_t \vect{q} + \matr{A}(\vect{q}) \partial_x \vect{q} = \SourceVector(\vect{q}),
\end{equation}
with hyperbolic transport matrix $\matr{A}(\vect{q}) \in \mathbb{R}^{(N+2)\times(N+2)}$ given by
\begin{equation} \label{e:hswmmat}
    \matr{A}(\vect{q}) =
    \begin{bmatrix}
        0 & 1 & &&& \\
        g h-u_m^2-\frac{1}{3}\alpha_1^2 & 2u_m & \frac{2}{3}\alpha_1 &&& \\
        -2u_m\alpha_1 & 2\alpha_1 & u_m & \frac{3}{5}\alpha_1 && \\
        -\frac{2}{3}\alpha_1^2 & 0 & \frac{1}{3}\alpha_1 & u_m & \ddots & \\
        &&& \ddots & \ddots & \frac{N+1}{2N+1}\alpha_1 \\
        &&&& \frac{N-1}{2N-1}\alpha_1 & u_m
    \end{bmatrix}.
\end{equation}
Note that the transport matrix $\matr{A}$ is a function that only depends on macroscopic $h$, $u_m$, plus $\alpha_1$, whereas it does not depend on the higher microscopic moments $\alpha_i$ for $i=2,\ldots,N$. This will be the first important ingredient allowing for efficient model reduction via POD-Galerkin in \cref{sec:POD-Galerkin} and DLRA in \cref{sec:DLRA} later. While the HSWME model appears simpler than SWME, it still retains the main nonlinearity of the SWME and is able to reproduce the complex nonlinear flow patterns of shallow flows, as demonstrated in numerous recent papers \cite{Garres2021,Huang2022,Pimentel2022,Amrita2022}.

The source term $\SourceVector (\vect{q}) \in \mathbb{R}^{N+2}$ of \cref{e:HSWME} is given by $\SourceVector (\vect{q}) = [0,{\sourceScalar}_0,\dots,{\sourceScalar}_{N}]^{T}$ and reads
\begin{equation}\label{e:source_term}
    {\sourceScalar}_i(\vect{q}) = -\frac{\nu}{\lambda} \left(2i+1\right) \left( u_m + \sum_{j=1}^N \alpha_j \right)  -\frac{\nu}{h} 4 \left(2i+1\right) \sum_{j=1}^N a_{i,j} \alpha_j,\quad i=0,\ldots,N,
\end{equation}
with slip length $\lambda$, viscosity $\nu$, and constants $a_{i,j}$ given by
\begin{equation}
    a_{i,j}=
    \begin{cases}
        0 \hspace{3cm} \quad\quad \textrm{ if } i+j =  \textrm{even},\\
        \frac{\min(i-1,j) \left( \min(i-1,j) +1\right)}{2} \quad \textrm{ if } i+j =  \textrm{odd}.
    \end{cases}
\end{equation}
Note that the first entry of the right-hand side friction term  $\SourceVector (\vect{q})$ is zero, leading to the conservation of mass, which is the integral of the water height $h$. More importantly, the other entries given by ${\sourceScalar}_i(\vect{q})$ only depend non-linearly on the water height $h$, while depending linearly on the remaining $u_m$ and $\alpha_i$. This will be another important ingredient for efficient model reduction to be exploited later.

\section{Space-time discretization with operator splitting}
\label{sec:discr}
The right-hand side friction term of \cref{e:HSWME} can be stiff for small $\lambda$ or small $h$. This prohibits the use of standard explicit schemes to solve the coupled model \cref{e:HSWME}. While stable explicit schemes like projective integration exist \cite{Amrita2022}, they typically require more time steps and some parameter tuning, which would be impractical for a robust model order reduction later. The left-hand side transport part is naturally discretised best with an explicit scheme using a standard CFL condition. We therefore split the full model \cref{e:HSWME} into a transport step and a friction step. We then treat the space-time discretization of both terms separately.

As explained in \cref{sec:Concept}, for a single time step of a numerical solution we apply a first order operator splitting to \cref{e:HSWME} in the same fashion as \cite{Huang2022}, in which we split the transport part \cref{e:transport} from the right-hand side friction part \cref{e:friction} as
\begin{alignat}{5}
      &\text{Step 1: transport} \qquad& \partial_t \vect{q}&+ \matr{A}(\vect{q}) \partial_x \vect{q} &\;=\;& 0\,, \label{e:transport2} \\
      &\text{Step 2: friction} \qquad& \partial_t \vect{q}& &\;=\;& \SourceVector (\vect{q})\,,\label{e:friction2}
\end{alignat}
where the solution of the transport and friction step will be considered separately in the following two subsections.

\subsection{A new macro-micro decomposition scheme for the transport step}
\label{sec:transport_step}
A straightforward discretization of the transport step \cref{e:transport} via a standard path-conservative numerical schemes \cite{Pimentel2022} is sufficient to solve the transport step of the full order model.
However, preparing for the model reduction in \cref{sec:POD-Galerkin} and \cref{sec:DLRA}, we here describe a new version exploiting the structure of the model using a macro-micro decomposition of the variables. This means additionally decomposing the vector of variables in the same way as they will be treated by the POD-Galerkin and DLRA methods in \cref{sec:POD-Galerkin} and \cref{sec:DLRA}, respectively.

\hspace{1em}
A standard numerical scheme as written in \cref{app1}, see \cref{e:transport_step_scheme}, considers an update of the full state vector. However, this does not allow to completely leverage the structure of the underlying model equations during the model reduction procedure later. To that end, we exploit the structure of the HSWME model \cref{e:HSWME}, by decomposing $\Q$ into two parts: (1) the first two (macroscopic) variables for the water height $h$ and the mean velocity $u_m$ (called $\U$) on the one hand and (2) the last $N$ (microscopic) variables for the moments $\alpha_i$ on the other hand (called $\V$). This can be written as
\begin{equation}
\label{eq-def:statematrixQU}
    \Q = \begin{bmatrix}
                    \U & \V
                 \end{bmatrix}, \quad
    \U = \begin{bmatrix}
                    h(x_j,t) & h(x_j,t) \, u_m(x_j,t)
                 \end{bmatrix}_{j} \in \mathbb{R}^{N_x \times 2},
\end{equation}
\begin{equation}
\label{eq-def:statematrixV}
    \V = \begin{bmatrix}
                    h(x_j,t)\, \alpha_1(x_j,t) & \hdots & h(x_j,t)\, \alpha_N(x_j,t)
                 \end{bmatrix}_{j} \in \mathbb{R}^{N_x \times N}.
\end{equation}

In the same fashion, we decompose the transport matrix $\matr{A}(\vect{q})$ into four blocks corresponding to the first two equations and the last $N$ equations and variables, respectively
\begin{equation}\label{e:block_matrices}
    \matr{A}(\vect{q}) =
        \begin{bmatrix}
           \matr{A}_{\vect{u}\vect{u}} & \matr{A}_{\vect{u}\vect{v}} \\ \matr{A}_{\vect{v}\vect{u}} & \matr{A}_{\vect{v}\vect{v}}
        \end{bmatrix},
\end{equation}
with blocks
\begin{alignat}{3}
    &\matr{A}_{\vect{u}\vect{u}} =
        \begin{bmatrix}
             & 1 \\
            g h-u_m^2-\frac{1}{3}\alpha_1^2 & 2u_m
        \end{bmatrix} \in \mathbb{R}^{2 \times 2},\quad
    &&\matr{A}_{\vect{u}\vect{v}} =
        \begin{bmatrix}
             & & \qquad\qquad\qquad\qquad\quad\,\,\,\\
            \frac{2}{3}\alpha_1 & & \quad
        \end{bmatrix} \in \mathbb{R}^{2 \times N},\\
    &\matr{A}_{\vect{v}\vect{u}} =
        \begin{bmatrix}
            \quad -2u_m\alpha_1 \quad & \quad 2\alpha_1 \\
            -\frac{2}{3}\alpha_1^2 & \\
            &                        \\
            &                        \\
            &
        \end{bmatrix} \in \mathbb{R}^{N \times 2},
    &&\matr{A}_{\vect{v}\vect{v}} =
        \begin{bmatrix}
            u_m & \frac{3}{5}\alpha_1 && \\
            \frac{1}{3}\alpha_1 & u_m & \ddots & \\
            & \ddots & \ddots & \frac{N+1}{2N+1}\alpha_1 \\
            && \frac{N-1}{2N-1}\alpha_1 & u_m
        \end{bmatrix} \in \mathbb{R}^{N \times N},
\end{alignat}
where all other entries are zeros.

We thus decompose the system into a set of macroscopic variables $U$ and another set of microscopic variables $V$. While the macroscopic variables are not reduced to allow structure preservation, the microscopic variables are used to generate an efficient reduced order model.

The space-time discretization of the transport part towards a solution scheme is described in the \cref{app1} for brevity. The derivation of reduced model by both POD-Galerkin as well as DLRA rely on the definition of this discretized scheme.  

After the solution of the transport step has been computed the scheme continues with the friction step, which is potentially stiff and requires a different, implicit scheme.

\subsection{Friction step}
\label{sec:friction_step}
To solve the space-homogeneous friction step \cref{e:friction}, an implicit scheme is necessary due to potential stiffness originating from small values for $\lambda$ or $h$. We will use the scheme from \cite{Huang2022}, but adopt it to our decomposition of variables $\Q$ into $\U$ and $\V$ as explained in this section.

The macro variables $\vect{u}=(h, h u_m)$ after the friction step are obtained using the definition of the source in \cref{e:source_term} and the two observations: (1) The height $h$ remains constant in time during the friction step. (2) The micro moments $\alpha_i$ remain constant during the macro update. This leads to the updated values $(h^{n+1}, h^{n+1} u_m^{n+1})$.

For the update of the remaining micro coefficients $\V$ during the micro friction step \cref{e:friction_micro}, we make use of the fact that $h$ and $u_m$ are now constant and $\V$ only occurs linearly. According to \cite{Huang2022}, \cref{e:friction} for the last $N$ moments can then be written as
\begin{align}\label{e:friction_ode}
    \dot{\V}_j &= \frac{1}{h_j^2}\SourceMatrix_1\V_j + \frac{1}{h_j}\SourceMatrix_2\V_j + u_{m,j}\SourceVector,
\end{align}
where $\SourceMatrix_1 \in \mathbb{R}^{N \times N}$ has entries
\begin{align}
    \SourceMatrix_{1,ij} = \begin{cases}
    -2\nu(2(i-1)+1)\cdot \min(i,j)\cdot (\min(i,j)+1) & \text{if } i+1+j \text{ is even}\\
    0 & \text{else}
    \end{cases}
\end{align}
and $\SourceMatrix_2= \SourceVector \cdot \mathbf{1}^{\top} \in \mathbb{R}^{N \times N}$, where $\SourceVector \in \mathbb{R}^{N}$ has entries ${\sourceScalar}_{i} = -\frac{\nu}{\lambda}(2i+1)$.

The micro friction step \cref{e:friction_ode} is then solved implicitly by a backward Euler method to overcome stability issues from potential stiffness. Note again that the water height $h_j$ is constant during the whole friction step and $u_{m,j}^{n+1}$ is given by the solution of the macro friction step \cref{e:friction_macro}. Applying the backward Euler method to \cref{e:friction_ode} yields
\begin{align}\label{e:friction_evolv}
    \V_j^{n+1} = \V_j^{n} + \Delta t \left( \frac{1}{h_j^2}\SourceMatrix_1\V_j^{n+1} + \frac{1}{h_j}\SourceMatrix_2\V_j^{n+1} + u_{m,j}^{n+1}\SourceVector \right).
\end{align}
We then define the matrix
\begin{align}
    \matr{D}_j = \matr{I}_{N} - \frac{\Delta t}{h_j^2}\SourceMatrix_1 - \frac{\Delta t}{h_j}\SourceMatrix_2
\end{align}
to arrive at the time update
\begin{align}\label{e:friction_Vupdate_step}
    \V_j^{n+1} = \matr{D}_j^{-1}\left(\V_j^{n} + \Delta t \, u_{m,j}^{n+1}\SourceVector\right).
\end{align}
We note that $\matr{D}_j^{-1}$ can be precomputed efficiently, as it does not depend on $u_m$ and $\alpha_{i}$. Thus, the updated micro variables in \cref{e:friction_Vupdate_step} can be computed efficiently without inverting a matrix during the online computation.

\section{Macro-micro decomposition for conservative POD-Galerkin}
\label{sec:POD-Galerkin}
The following section addresses the POD-Galerkin reduction of the HSWME model from \cref{sec:HSWME}.
As already pointed out in \cref{sec:Concept} the presented approach applies the reduction only to the microscopic higher-order moments $\vect{v}=[h\alpha_1,\dots,h\alpha_N]^\top\in\mathbb{R}^N$.
Furthermore, in contrast to the conventional POD-Galerkin approach that reduces the full state-space of the discretized PDE, we only reduce the dimensions of the PDEs' state-components. This leads to a decoupling between the time-space dynamics and the correlation between the components and avoids the separation of spatial-temporal dynamics that can lead to slow decaying approximation errors, known as the Kolmogorov $N$-width problem \cite{OhlbergerRave2015,GreifUrban2019}. 
For a typical dam break scenario we showcase this slow decay in \cref{fig:water-column:error-decay-pod} of our numerical \cref{sec:numerics}.
To obtain a reduced model that allows rapid yet accurate predictions over a range of different parameters, one often uses a two phase offline-online procedure, which is explained in the following.

\subsection{Offline phase}
\label{subsec:offlinephase}
In the \emph{offline phase} the reduced basis and the operators of the reduced model space are precomputed.
Note that in contrast to the DLRA approach from \cref{sec:DLRA} one is willing to forego large up-front offline costs in favor of a more efficient reduced model for the online phase.
The reduced basis is formed by collecting snapshots of the solution of \cref{e:HSWME} at $N_\mathrm{t}$ different time/parameter instances
\begin{equation}
    \mathcal{V}=\{\vect{v}(x,t_1),\dots,\vect{v}(x,t_{\Nt-1})\}\,.
\end{equation}
The POD approximates the microscopic higher moment vector using an orthonormal basis $\{\vect{w}\}_{k=1,\dots,r}$
\begin{align}
    \vect{v}(x,t)\approx \widetilde{\vect{v}}(x,t)&=\sum_{k=1}^r \hat{\alpha}_k(x,t) \vect{w}_k\,, \qquad r\ll N\le \Nt\\
                                             &=\matr{W}\hat{\vect{v}}(x,t)\,,
\end{align}
where $\matr{W}=[\vect{w}_1,\dots,\vect{w}_r]\in\mathbb{R}^{N\times r}$ collects the basis vectors and $\hat{\vect{v}}=[\hat{\alpha}_1,\dots,\hat{\alpha}_r]^\top\in\mathbb{R}^r$ denotes the reduced space-time coefficients.
In a fully discrete setting this basis is computed by a truncated \emph{singular value decomposition} (SVD) of the \emph{snapshot matrix}
\begin{align}
    \matr{V}^\mathrm{POD}
        =\begin{bmatrix}
        \matr{V}^0\\
        \vdots\\
        \matr{V}^{\Nt-1}
        \end{bmatrix}\in\mathbb{R}^{(N_x\Nt)\times N}\,,
\end{align}
where the $\matr{V}^n$ are defined in \cref{eq-def:statematrixV}. The truncated SVD of $\matr{V}^\mathrm{POD}$ yields
\begin{align}
\label{eq:svd-POD}
   \matr{V}^\mathrm{POD}\approx \V^\mathrm{POD}_r = \matr{\Psi\Sigma}\matr{W}^\top
\end{align}
  Here, $\matr{\Sigma}=\diag{(\sigma_1,\dots,\sigma_r)},\,r\ll M$, is a diagonal matrix containing the largest $r$ singular values $\sigma_1\ge\sigma_2\ge\dots \ge \sigma_r$ and $\matr{\Psi}\in\mathbb{R}^{(N_x\Nt)\times r}$, $\matr{W}\in\mathbb{R}^{N\times r}$ are orthogonal matrices containing the left and right singular vectors, respectively. The latter are also termed \emph{modes} in the following.
  
According to the Eckart-Young-Mirsky theorem \cite{EcY1936,Mirsky1960} $\V_r^\mathrm{POD}$ is the best rank $r$ approximation and the resulting error in the Frobenius norm is rigorously computed from the trailing singular values
\begin{equation}
    \norm{\V^\mathrm{POD}-\V_r^\mathrm{POD}}_\mathrm{F}^2=\sum_{k=r+1}^m\sigma_k^2\,.
\end{equation}
A common choice for $r$ is to truncate after a certain energy percentage (e.g. $E_\mathrm{cum}\ge 95\%$) is reached in the reduced system compared to the full system:
\begin{equation}\label{e:ePOD}
    E_\mathrm{cum}=\frac{\norm{\V^\mathrm{POD}_r}_\mathrm{F}}{\norm{\V^\mathrm{POD}}_\mathrm{F}}=\frac{\sum_{k=1}^r\sigma_k^2}{\sum_{k=1}^m\sigma_k^2}.
\end{equation}
  
Calculating an SVD can be challenging, particularly when dealing with higher-dimensional space, due to the potentially prohibitive size of the snapshot matrix. This can be circumvented by either using randomized or wavelet techniques to calculate the POD-modes \cite{KrahEngelsSchneiderReiss2022,HalkoMartinssonTropp2011} more efficiently or applying POD-greedy sampling methods \cite{HaasdonkOhlberger2008,WaldherrHaasonk2012} that reduce the number of snapshots needed to form a POD basis.
  
In the context of Galerkin projections the column space of $\matr{W}$ yields the trial space
\begin{equation}
     \V(t)\approx \widetilde{\matr{V}}(t) = \matr{W}\Vhat(t), \qquad \Vhat\colon [0,T]\to \mathbb{R}^{r\times N_x}\,.
\end{equation}

\subsection{Online Phase}
    As explained in \cref{sec:Concept}, the \emph{online phase} evolves the dynamics in two steps. In the first step we evolve the first two macro variables $\vect{U}$ without any model order reduction and in the second step we use the evolved $\vect{U}$ as an input of our reduced system of the microscopic higher order moments $\V$. The reduced system is gathered by projecting
    \cref{eq:RHSV-spacediscrete,e:friction_evolv} onto the test space which is the reduced subspace spanned by the columns of $\matr{W}$.

\subsubsection{Transport step}
    For the transport part in \cref{eq:RHSV-spacediscrete} we define the reduced transport term as
    \begin{equation}
        \widehat{\advectV}(\Vhat):=\matr{W}^\top\advectV(\U,\matr{W}\Vhat(t))\in\mathbb{R}^{r\times N_x}\,.
    \end{equation}
     Evaluating this at the $j$th cell we obtain:
    \begin{equation}
        \label{e:GalerkinROM_advect}
    \begin{aligned}
         \relax [\widehat{\advectV}(\Vhat)]_j &=\frac{1}{2\Delta t}(\Vhat_{j-1}-2\Vhat_{j}+\Vhat_{j+1})\\
         &- \frac{1}{2\Delta x} \left( \hat{\matr{A}}_{j+1,\vect{v}\vect{u}} \left(\U_{j+1} - \U_{j}\right) + \hat{\matr{A}}_{j+1,\vect{v}\vect{v}} \left(\Vhat_{j+1} - \Vhat_{j}\right)\right) \\
        &- \frac{1}{2\Delta x} \left( \hat{\matr{A}}_{j,\vect{v}\vect{u}} \left(\U_{j} - \U_{j-1}\right) + \hat{\matr{A}}_{j,\vect{v}\vect{v}} \left(\Vhat_{j} - \Vhat_{j-1}\right) \right),
    \end{aligned}
    \end{equation}
    where we define the reduced non-linear operators  $\hat{\matr{A}}_{j,\vect{v}\vect{u}} := \mathbf{W}^{\top}\matr{A}_{j,\vect{v}\vect{u}}$ and $\hat{\matr{A}}_{j,\vect{v}\vect{v}}$, where:
    \begin{align}
    \label{e:GalerkinROM_reducedOP_A}
          \hat{\matr{A}}_{j,\vect{v}\vect{v}} &:= \mathbf{W}^{\top} \matr{A}_{j,\vect{v}\vect{v}}\mathbf{W}=\left(\frac12(\alpha_{1,j+1}+\alpha_{1,j})\mathbf{\hat A}+\frac12(u_{m,j+1}+u_{m,j})\mathbf{I}_r\right)\\
          \mathbf{\hat A}&:= \mathbf{W}^{\top}\mathbf{A}\mathbf{W}\in\mathbb{R}^{r\times r},
    \end{align}
    with identity matrix $\matr{I}_r \in \mathbb{R}^{r \times r}$.
        Note, since $\mathbf{\hat A}$ is precomputed in the offline phase, the non-linear term \cref{e:GalerkinROM_reducedOP_A} only has to be evaluated for $r$ instead of $N$ components, which generates the speedup.
        To further simplify the non-linear terms POD is often used in combination with sparse sampling
methods that sample the nonlinear terms at a few components to approximate them in a low-dimensional space. Examples are the \emph{discrete empirical interpolation method} (DEIM) \cite{ChaturantabutSorensen2010}, that relates back to the empirical interpolation method (EIM) \cite{barrault2004empirical,grepl2007efficient}, the gappy POD \cite{everson1995karhunen,bui2003proper,AstridWeilandWillcoxBackx2008}, or the \emph{energy-conserving sampling and weighting} (ECSW) method \cite{farhat2014dimensional}.
However, we refrain from using them here for a direct comparison with DLRA, that does not make use of sparse sampling methods, yet.

\subsubsection{Friction step}
    Similar to the transport term we evaluate the micro friction term \cref{e:friction_micro} at the $j$th component projected onto the test space. We define the components of the reduced source term $\hat{\SourceMatrix}(\hat{\matr{V}}):=\matr{W}^\top{\SourceMatrix}(\matr{U},\matr{W}\hat{\matr{V}})\in\mathbb{R}^{r\times N_x}$ by
    \begin{align}
        [\hat{\SourceMatrix}(\hat{\matr{V}})]_j &= \frac{1}{h_j^2} \hat{\SourceMatrix}_1 \hat{\matr{V}}_j + \hat{\SourceMatrix}_2 \hat{\matr{V}}_j + u_{m,j}\hat{\SourceVector}\,\\
        \hat{\SourceMatrix}_1 &= \matr{W}^\top \SourceMatrix_1 \matr{W}\in\mathbb{R}^{r\times r}, \, \qquad \hat{\SourceMatrix}_2 = \matr{W}^\top \SourceMatrix_2 \matr{W} \in\mathbb{R}^{r\times r}, \qquad \hat{\SourceVector}=\matr{W}^\top \SourceVector \in\mathbb{R}^r
    \end{align}
    The complexity of the source term computation is reduced, since only $r$ components have to be evaluated.

\section{Dynamical low-rank approximation}
\label{sec:DLRA}
As an alternative to POD-Galerkin, we propose to evolve the microscopic higher-order moments with DLRA introduced in \cite{koch2007dynamical}. This method is data-driven in the sense that it closes the moment equations without assuming physical properties, based on the real-time solution data. 

\subsection{Macro-micro decomposition for dynamical low-rank approximation}

The core idea of DLRA is to evolve the solution on a low-rank manifold. That is, DLRA represents the solution as a low-rank factorization and provides evolution equations for the individual factors. Therefore, DLRA can be interpreted as a Galerkin method which updates not only the expansion coefficients but also the basis functions in time. To preserve the structure of water height and momentum, we apply DLRA to the microscopic correction terms $v_i := h\alpha_i$ only. Collecting the spatially discretized correction terms in a matrix $\V(t)\in\mathbb{R}^{N_x\times N}$, where $v_{ji} = h(t,x_j)\alpha_i(t,x_j)$, we define a low-rank approximation as $\V(t) = \mathbf{X}(t)\mathbf{S}(t)\mathbf{W}(t)^{\top}$, where $\mathbf{X}\in\mathbb{R}^{N_x\times r}$ and $\mathbf{W}\in\mathbb{R}^{N\times r}$ can be interpreted as the collection of $r$ basis vectors in space and moment order with a corresponding coefficient matrix $\mathbf{S}\in\mathbb{R}^{r\times r}$.

To benefit from this representation, we want to ensure that the method works on these factors only and never needs to compute and store the full solution $\V$. To preserve the low-rank structure of the solution, we therefore force the solution at all times $t$ to remain in the manifold of rank $r$ matrices, which we call $\mathcal{M}_r$. This can be ensured when the time derivative of $\V$ lies in the tangent space of rank $r$ matrices, i.\,e., the solution when advancing in time does not leave the manifold $\mathcal{M}_r$. We denote the tangent space at $\V(t)$ as $\mathcal{T}_{\V(t)}\mathcal{M}_r$. Then, the time evolution equations for the basis vectors and coefficients must satisfy
\begin{align}\label{eq:DLRAcond}
    \dot{\V}(t) \in\mathcal{T}_{\V(t)}\mathcal{M}_r \qquad \text{ such that}\quad \Vert \dot{\V}(t) - \RhsV(\U(t),\V(t)) \Vert \rightarrow \text{min!}
\end{align}
The first condition conserves the representation $\V(t) = \mathbf{X}(t)\mathbf{S}(t)\mathbf{W}(t)^{\top}$ and the second condition minimizes the residual, which is essentially a Galerkin condition. Hence, an equivalent formulation of \eqref{eq:DLRAcond} is
\begin{align}
    \left\langle \dot{\V}(t) - \RhsV(\U(t),\V(t)) , \delta \V\right\rangle = 0 \quad \forall \delta \V \in\mathcal{T}_{\V(t)}\mathcal{M}_r.
\end{align}
This represents a Galerkin method, which chooses test functions based on the solution data and the geometry of the manifold of rank $r$ matrices. An admissible choice of $\delta \V$ when the solution reads $\V(t) = \mathbf{X}(t)\mathbf{S}(t)\mathbf{W}(t)^{\top}$ is $\delta \V = \mathbf{X}_i\mathbf{W}_j^{\top}$ which indeed lies in the tangent space. Here, $\mathbf{X}_i$ and $\mathbf{W}_j$ denote the $i_{\mathrm{th}}$ and $j_{\mathrm{th}}$ columns of the basis matrices, respectively. This test function yields an evolution equation for the coefficient matrix
\begin{align}
    \dot{\mathbf{S}}(t) = \mathbf{X}(t)^{\top}\RhsV(\mathbf{X}(t)\mathbf{S}(t)\mathbf{W}(t)^{\top})\mathbf{W}(t).
\end{align}
Choosing test functions $\delta \V = \mathbf{X}_i$ and $\delta \V = \mathbf{W}_i$, which again lie in the tangent space, results in evolution equations for the respective basis matrices
\begin{align}
    \dot{\mathbf{X}}(t) =& \,(\mathbf{I} - \mathbf{X}(t)\mathbf{X}(t)^{\top})\RhsV(\mathbf{X}(t)\mathbf{S}(t)\mathbf{W}(t)^{\top})\mathbf{W}(t)\mathbf{S}(t)^{-1},\\
    \dot{\mathbf{W}}(t) =& \,(\mathbf{I} - \mathbf{W}(t)\mathbf{W}(t)^{\top})\RhsV(\mathbf{X}(t)\mathbf{S}(t)\mathbf{W}(t)^{\top})^{\top}\mathbf{X}(t)\mathbf{S}(t)^{-\top}.
\end{align}
Solving these evolution equations with classical time integration schemes is inefficient, since the coefficient matrices are often ill conditioned and the time step size is dictated by the smallest absolute eigenvalue of $\mathbf{S}$. Two robust time integrators, which guarantee stability, are the \emph{projector--splitting integrator} \cite{lubich2014projector} as well as the BUG integrator \cite{ceruti2020time}. Their main strategy is to not evolve $\mathbf{X}$ and $\mathbf{W}$ in time directly, but to evolve a linear transformation $\mathbf{K} := \mathbf{X}\mathbf{S}$ and $\mathbf{L} := \mathbf{W}\mathbf{S}^{\top}$ and retrieve the basis matrices through a QR-decomposition. In this work, we focus on the BUG integrator, which consists of three update steps
\begin{enumerate}
    \item \textbf{$K$-step}: Update $\mathbf X^{0}$ to $\mathbf X^{1}$ via
    \begin{linenomath*}\begin{align}
        \dot{\mathbf K}(t) &= \RhsV(\mathbf{K}(t)\mathbf{W}^{0,\top})\mathbf{W}^0\;, \qquad \mathbf K(t_0) = \mathbf{X}^0\mathbf{S}^0\;.\label{eq:KStepSemiDiscreteUI}
    \end{align}\end{linenomath*}
Determine $\mathbf X^1$ with $\mathbf K(t_1) = \mathbf X^1 \mathbf R$ and store $\mathbf M = \mathbf X^{1,\top}\mathbf X^0$.
\item \textbf{$L$-step}: Update $\mathbf W^0$ to $\mathbf W^1$ via
\begin{linenomath*}\begin{align}
\dot{\mathbf L}(t) &= \RhsV(\mathbf{X}^0\mathbf{L}(t)^T)^T\mathbf{X}^{0}\;, \qquad \mathbf L(t_0) = \mathbf{W}^{0}\mathbf{S}^{\top}\;.\label{eq:LStepSemiDiscreteUI}
\end{align}\end{linenomath*}
Determine $\mathbf W^1$ with $\mathbf L(t_1) = \mathbf W^1\mathbf{\widetilde R}$ and store $\mathbf N = \mathbf W^{1,\top} \mathbf W^0$.
\item \textbf{$S$-step}: Update $\mathbf S^0$ to $\mathbf S^1$ via
\begin{linenomath*}\begin{align}
\dot{\mathbf S}(t) = \mathbf{X}^{1,\top}\RhsV(\mathbf{X}^{1}\mathbf{S}(t)\mathbf{W}^{1,\top})\mathbf{W}^1\;, \qquad \mathbf S(t_0) &= \mathbf M\mathbf S^0 \mathbf N^{\top}\label{eq:SStepSemiDiscreteUI}
\end{align}\end{linenomath*}
and set $\mathbf{S}^1 = \mathbf S(t_1)$.
\end{enumerate}
The solution at the next time step is then given by $\V(t_1) = \mathbf{X}^1\mathbf{S}^1\mathbf{W}^{1,\top}$. Note that the first two equations which evolve the basis in time can be updated in parallel followed by a serial update of the coefficient vector.
Effectively the size of the original matrix ODE is reduced from a $N\times N_x$ system to three smaller matrix ODEs of size: $N_x\times r$ ($K$-step), $N\times r$ ($L$-step), and $r\times r$ ($S$-step). We wish to underline the requirements of a user-determined choice of the rank $r$, which needs to be determined such that the solution exhibits a sufficiently accurate approximation at minimal computational costs and memory requirements. Note however that the BUG integrator allows for an extension to rank adaptivity \cite{ceruti2022rank}. Here, the basis after the $K$ and $L$ steps is augmented with the basis at time $t_0$ to enlarge the approximation space from $r$ to $2r$ basis vectors. After a $2r\times 2r$ coefficient update in the $S$-step, the solution is truncated to a new rank $r_1 \leq 2r$ according to a user-determined tolerance parameter $\vartheta$. While this parameter needs to be chosen before the simulation, it has a clear interpretation as the truncation error in every time step. For clarity of presentation, we focus on the fixed-rank BUG integrator in our derivations, the numerical results however include both fixed-rank and rank-adaptive computations. For more information on rank-adaptivity and further properties of the rank-adaptive integrator, we refer to \cite{ceruti2022rank}. Further approaches for rank-adaptivity in dynamical low-rank approximation and model order reduction are, for example, \cite{hochbruck2023rank,ceruti2023parallel,dektor2021rank,hesthaven2022rank,ehrlacher2017dynamical}.

The classical dynamical low-rank approximation approach often does not preserve important physical properties. This stems from the fact that DLRA can remove basis functions which are needed for conservation. However, problem-dependent adaptations to the classical DLRA integrators can provide conservation properties. As an example, conservation of solution invariants up to a tolerance parameter can be achieved with a basis augmentation step \cite{ceruti2022rank}. Moreover, \cite{einkemmer2022robust} uses a basis augmentation as well as a reformulation of the $K$, $L$ and $S$-step to preserve mass, momentum and energy in the Vlasov equations. Our approach in this work ensures local mass conservation by decomposing the dynamics of the macroscopic (conserved) water height and momentum from the dynamics of the microscopic correction terms. This means that the conservation law structure of the water height equation is not altered by DLRA.

\subsection{Evolution equations for the low-rank HSWME}
In the following we derive efficient representations of the $K$, $L$, and $S$-steps for the micro transport step \cref{e:transport_micro} and the micro friction step \cref{e:friction_micro} of the HSWME model. Note that the macroscopic variables are not altered by the DLRA and computed during the macro transport and macro friction step as explained in \cref{sec:discr}. Note that the derivation of evolution equations is performed on the semi-discrete level, i.e., the system is already discretized in the vertical direction $z$ by a moment approximation and the horizontal direction by a finite volume discretization. DLRA can, however, also be derived in the fully continuous formulation \cite{einkemmer2019quasi}, leading to a continuous set of low-rank evolution equations. When using BUG integrators combined with an equivalent moment and finite volume discretization of the continuous low-rank evolution equations, these two approaches are equivalent, and we choose the discrete formulation for ease of presentation.

\subsubsection{Transport step}
For the micro transport part \cref{e:transport_micro}, the spatially discretized right-hand side \cref{eq:RHSV-spacediscrete} is given by
\begin{equation}
\begin{aligned}
    \relax [\advectV(\U,\V)]_j &=\frac1{2\Delta t}(\V_{j-1}-2\V_{j}+\V_{j+1})\\
    &- \frac{1}{2\Delta x} \left( \matr{A}_{j+1,\vect{v}\vect{u}} \left(\U_{j+1} - \U_{j}\right) + \matr{A}_{j+1,\vect{v}\vect{v}} \left(\V_{j+1} - \V_{j}\right)\right) \\
    &- \frac{1}{2\Delta x} \left( \matr{A}_{j,\vect{v}\vect{u}} \left(\U_{j} - \U_{j-1}\right) + \matr{A}_{j,\vect{v}\vect{v}} \left(\V_{j} - \V_{j-1}\right) \right).
\end{aligned}
\end{equation}
\textbf{$K$-step:} To derive the $K$-step, we denote the $j$-th row of $\mathbf{K}$ as $\mathbf{K}_j(t)\in\mathbb{R}^{r}$ and represent the solution $\V$ at spatial cell $j$ as $\V_j(t) = \mathbf{W}^0\mathbf{K}_j(t)$. Multiplying the right-hand side $\advectV$ with $\mathbf{W}^{0,\top}$ then yields
\begin{equation}
\begin{aligned}
    \dot{\matr{K}}_j &=\frac1{2\Delta t}(\matr{K}_{j-1}-2\matr{K}_{j}+\mathbf{K}_{j+1}) \\
    &- \frac{1}{2\Delta x} \left( \matr{W}^{0,\top}\matr{A}_{j+1,\vect{v}\vect{u}} \left(\U_{j+1} - \U_{j}\right) + \widetilde{\matr{A}}_{j+1,\vect{v}\vect{v}} \left(\matr{K}_{j+1} - \matr{K}_{j}\right)\right) \\
    &- \frac{1}{2\Delta x} \left( \matr{W}^{0,\top}\matr{A}_{j,\vect{v}\vect{u}} \left(\U_{j} - \U_{j-1}\right) + \widetilde{\matr{A}}_{j,\vect{v}\vect{v}} \left(\matr{K}_{j} - \matr{K}_{j-1}\right) \right),
\end{aligned}
\end{equation}
where $\widetilde{\matr{A}}_{j,\vect{v}\vect{v}} := \matr{W}^{0,\top} \matr{A}_{j,\vect{v}\vect{v}}  \matr{W}^{0}$. Note that since $\widetilde{\matr{A}}_{j,\vect{v}\vect{v}} \in \mathbb{R}^{r\times r}$ and $\matr{W}^{0,\top}\matr{A}_{j+1,\vect{v}\vect{u}} \in \mathbb{R}^{r\times 2}$, the main memory requirements stem from storing $\matr{K}_j$ at all spatial cells $j$, i.\,e., memory requirements are $\mathcal{O}(r\cdot N_x)$ opposed to the original method's requirements of $\mathcal{O}(N_x\cdot N)$. An efficient computation of $\widetilde{\matr{A}}_{j,\vect{v}\vect{v}}$ precomputes $\mathbf{\widetilde A}:= \mathbf{W}^{0,\top}\mathbf{A}\mathbf{W}^0$ which requires $\mathcal{O}(r^2 \cdot N)$ operations. Using the matrix $\mathbf{\widetilde A}\in\mathbb{R}^{r\times r}$ to compute $\widetilde{\matr{A}}_{j,\vect{v}\vect{v}}$ via
\begin{align}
    \widetilde{\matr{A}}_{j,\vect{v}\vect{v}} = \left(\frac12(\alpha_{1,j}^{}+\alpha_{1,j-1}^{})\mathbf{\widetilde A}+\frac12(u_{m,j}^{}+u_{m,j-1}^{})\mathbf{I}_r\right)
\end{align}
requires $O(r^2 \cdot N_x)$ operations. Hence, the computational costs for the $K$-step are $C_K \lesssim r^2\cdot(N_x+N)$.\\

\textbf{$L$-step:} To derive the $L$-step, we represent the solution $\V$ at spatial cell $j$ as $\V_j(t) = \mathbf{L}(t)\mathbf{X}_{j}^0$. Moreover, we test with $\mathbf{X}^0$, meaning that we multiply the right-hand side with $\mathbf{X}^0_{j}$ and sum over $j$. Hence, using Einstein's sum convention and writing out $\widetilde{\matr{A}}_{j,\vect{v}\vect{v}}$ we have
\begin{equation}
\begin{aligned}
    \dot{\mathbf{L}} =& \frac1{2\Delta t}\mathbf{L}(t)\cdot(\mathbf{X}_{j-1}^0-2\mathbf{X}_{j}^0+\mathbf{X}_{j+1}^0)\mathbf{X}_{j}^{0,\top}\\
    &- \frac{1}{4\Delta x} \left(\mathbf{A}\mathbf{L}(t) \cdot(\alpha_{1,j+1}+\alpha_{1,j})+\mathbf{L}(t)\cdot(u_{m,j+1}+u_{m,j})\right)\cdot(\mathbf{X}^0_{j+1}-\mathbf{X}^0_{j})\mathbf{X}_{j}^{0,\top}\\
    &- \frac{1}{4\Delta x} \left(\mathbf{A}\mathbf{L}(t) \cdot(\alpha_{1,j}+\alpha_{1,j-1})+\mathbf{L}(t)\cdot(u_{m,j}+u_{m,j-1})\right)\cdot(\mathbf{X}^0_{j}-\mathbf{X}^0_{j-1})\mathbf{X}_{j}^{0,\top}\\
     &- \frac{1}{2\Delta x} \left( \matr{A}_{j+1,\vect{v}\vect{u}} (\U_{j+1}-\U_{j})\mathbf{X}_{j}^{0,\top}+ \matr{A}_{j,\vect{v}\vect{u}} ( \U_{j}-\U_{j-1})\mathbf{X}_{j}^{0,\top} \right)\,.
\end{aligned}
\end{equation}
To simplify the structure of the resulting equations, we define
\begin{align}
    \mathbf{\widetilde X}^0 :=& \sum_j(\mathbf{X}_{j-1}^0-2\mathbf{X}_{j}^0+\mathbf{X}_{j+1}^0)\mathbf{X}_{j}^{0,\top},\\
    \bm{\widetilde \alpha}_+^0 :=& \sum_j(\alpha_{1,j+1}+\alpha_{1,j})(\mathbf{X}^0_{j+1}-\mathbf{X}^0_{j})\mathbf{X}_{j}^{0,\top},\\
    \bm{\widetilde u}_{+}^0 :=& \sum_j(u_{m,j+1}+u_{m,j})(\mathbf{X}^0_{j+1}-\mathbf{X}^0_{j})\mathbf{X}_{j}^{0,\top},
\end{align}
and $\bm{\widetilde \alpha}_-^0$ as well as $\bm{\widetilde u}_-^0$ accordingly. This results in
\begin{equation}
\begin{aligned}
    \dot{\mathbf{L}} =& \frac1{2\Delta t}\mathbf{L}(t)\mathbf{\widetilde X}^0- \frac{1}{4\Delta x} \left(\mathbf{A}\mathbf{L}(t) \left(\bm{\widetilde \alpha}_+^0+\bm{\widetilde \alpha}_-^0\right)+\mathbf{L}(t)\left(\bm{\widetilde u}_+^0+\bm{\widetilde u}_-^0\right)\right)\\
     &- \frac{1}{2\Delta x} \left( \matr{A}_{j+1,\vect{v}\vect{u}} (\U_{j+1}-\U_{j})\mathbf{X}_{j}^{0,\top}+ \matr{A}_{j,\vect{v}\vect{u}} ( \U_{j}-\U_{j-1})\mathbf{X}_{j}^{0,\top} \right).
\end{aligned}
\end{equation}
Note that the highest memory requirements come from storing $\mathbf{L}\in\mathbb{R}^{N\times r}$. Computing $\mathbf{\widetilde X}^0$, $\bm{\widetilde \alpha}_{\pm}^0$ and $\bm{\widetilde u}_{\pm}^0$ requires $\mathcal{O}(r^2 \cdot N_x)$ operations and computing $\mathbf{A}\cdot \mathbf{L}$ requires $\mathcal{O}(N\cdot r)$ operations, since $\mathbf{A}$ only has off-diagonal entries. The multiplication of $\mathbf{L}$ and $\mathbf{\widetilde X}^0$, $\bm{\widetilde \alpha}_{\pm}^0$ as well as $\bm{\widetilde u}_{\pm}^0$ requires $\mathcal{O}(r^2\cdot N)$ operations, i.\,e., we have a computational cost for the $L$-step of $C_L \lesssim O(r^2 \cdot (N_x + N))$.\\

\textbf{$S$-step:} To derive the $S$-step, we represent the microscopic solution $\V$ at spatial cell $j$ as $\V_j(t) = \mathbf{W}^1\mathbf{S}(t)^{\top}\mathbf{X}_{j}^1$. Moreover, we test with $\mathbf{X}^1$ and $\mathbf{W}^1$, i.\,e., we multiply the right-hand side with $\mathbf{X}^1_{j}$ and sum over $j$ as well as multiply with $\mathbf{W}^{1,\top}$. With the previous definitions we have
\begin{equation}
\begin{aligned}
    \dot{\mathbf{S}} =& \frac1{2\Delta t}\mathbf{S}(t)\mathbf{\widetilde X}^1- \frac{1}{4\Delta x} \left(\mathbf{W}^{1,\top}\mathbf{A}\mathbf{W}^{1}\mathbf{S}(t)^{\top} \left(\bm{\widetilde \alpha}_+^1+\bm{\widetilde \alpha}_-^1\right)+\mathbf{S}(t)^{\top}\left(\bm{\widetilde u}_+^1+\bm{\widetilde u}_-^1\right)\right)\\
     &- \frac{1}{2\Delta x} \left( \mathbf{W}^{1,\top} \matr{A}_{j+1,\vect{v}\vect{u}} (\U_{j+1}-\U_{j})\mathbf{X}_{j}^{1,\top}+\mathbf{W}^{1,\top} \matr{A}_{j,\vect{v}\vect{u}} ( \U_{j}-\U_{j-1})\mathbf{X}_{j}^{1,\top} \right).
\end{aligned}
\end{equation}
The memory requirements are $\mathcal{O}(r^2)$ and following the discussion for the $K$ and $L$-steps, the computational costs for the $S$-step again are $C_S \lesssim \mathcal{O}(r^2 \cdot (N_x + N))$.

\subsubsection{Friction step}
For the friction step \cref{e:friction_micro} we have according to \cref{e:friction_ode}
\begin{align}
    \dot{\V}_j = \frac{1}{h_j^2} \SourceMatrix_1 \V_j+\frac{1}{h_j} \SourceMatrix_2 \V_j + u_{m,j} \SourceVector.
\end{align}
For simplicity of notation, we define the matrix \emph{ordinary differential equation} (ODE)
\begin{align}\label{eq:frictionFull}
    \dot{\V} = \mathbf{h}^{-2}\V \SourceMatrix_1^{\top}+\mathbf{h}^{-1}\V \SourceMatrix_2^{\top} + \mathbf{u}_{m} \SourceVector^{\top},
\end{align}
where $\mathbf{h} = \text{diag}(h_1,\cdots,h_{N_x})$ and $\mathbf{u}_m = \text{diag}(u_{m,1},\cdots,u_{m,N_x})$ will be the values obtained from a previous macro friction update \cref{e:friction_macro} and assumed known during the micro step \cref{e:friction_micro}. Let us again derive $K$, $L$ and $S$-steps and directly define a time discretization.\\

\textbf{$K$-step:} Let us use the representation $\V = \mathbf{K}(t)\mathbf{W}^{1,\top}$ and test with $\mathbf{W}^1$. Then, the $K$-step equation reads
\begin{align}
     \dot{\mathbf{K}} = \mathbf{h}^{-2}\mathbf{K}\mathbf{W}^{1,\top} \SourceMatrix_1^{\top}\mathbf{W}^1+\mathbf{h}^{-1}\mathbf{K}\mathbf{W}^{1,\top} \SourceMatrix_2^{\top}\mathbf{W}^1 + \mathbf{u}_{m}^{} \SourceVector^{\top}\mathbf{W}^1.
\end{align}
Since friction terms are commonly stiff, we use an implicit Euler method to discretize in time. That is, for a fixed spatial cell $j$ and defining $\widehat{\SourceMatrix}_1 := \mathbf{W}^{1,\top}\SourceMatrix_1\mathbf{W}^1$ as well as $\widehat{\SourceMatrix}_2 := \mathbf{W}^{1,\top}\SourceMatrix_2\mathbf{W}^1$ we have
\begin{align}
     \left(\mathbf I_r -\Delta th_j^{-2} \widehat{\SourceMatrix}_1 -\Delta th_j^{-1} \widehat{\SourceMatrix}_2 \right)\mathbf{K}_j^{n+1} = \mathbf{K}_j^{n} + \Delta tu_{m,j}^{n+1} \SourceVector^{\top}\mathbf{W}^1,
\end{align}
where $u_{m,j}^{n+1}$ is the updated velocity from the macro friction step \cref{e:friction_macro}.

We again retrieve the time updated spatial basis which we denote by $\mathbf{X}^2$ by a QR decomposition of $\mathbf{K}^{n+1}$. Hence, at $N_x$ spatial cells, we need to solve a linear system with $r$ unknowns. Moreover, we have to compute flux matrices with $r^2$ entries which require $N$ operations per entry. Thus, the computational costs are $C_K^f \lesssim N_x \cdot r^3 + N\cdot r^2$.\\

\textbf{$L$-step:} Use the representation $\V = \mathbf{X}^1\mathbf{L}^{\top}$ and test the transposed of \eqref{eq:frictionFull} with $\mathbf{X}^1$. Then, the $L$-step equation reads
\begin{align}
    \dot{\mathbf{L}} = \SourceMatrix_1 \mathbf{L}\mathbf{X}^{1,\top}\mathbf{h}^{-2}\mathbf{X}^1+ \SourceMatrix_2 \mathbf{L}\mathbf{X}^{1,\top}\mathbf{h}^{-1}\mathbf{X}^1 + \SourceVector \mathbf{u}_{m}^{\top}\mathbf{X}^1\,.
\end{align}
Let us define $\mathbf{\widehat h}_1^{-2}:=\mathbf{X}^{1,\top}\mathbf{h}^{-2}\mathbf{X}^1$ and $\mathbf{\widehat h}_1^{-1}:=\mathbf{X}^{1,\top}\mathbf{h}^{-1}\mathbf{X}^1$ and again use an implicit Euler time discretization. This gives
\begin{align}
    \mathbf{L}^{n+1} -\Delta t \SourceMatrix_1 \mathbf{L}^{n+1}\mathbf{\widehat h}_1^{-2}-\Delta t \SourceMatrix_2 \mathbf{L}^{n+1}\mathbf{\widehat h}_1^{-1}=  \mathbf{L}^{n}+\Delta t \SourceVector \mathbf{u}_{m}^{n+1,\top}\mathbf{X}^1.
\end{align}
The time updated coefficient basis which we denote by $\mathbf{W}^2$ is then obtained by a QR decomposition of $\mathbf{L}^{n+1}$. Hence, we need to solve an $r\cdot N$ system of linear equations which requires $\mathcal{O}(r^3\cdot N^3)$ operations. Moreover, computing $\mathbf{\widehat h}_1^{-1}$ and $\mathbf{\widehat h}_1^{-2}$ requires $\mathcal{O}(N_x\cdot r^2)$ operations, hence $C_L^{f}\lesssim r^3\cdot N^3+N_x\cdot r^2$.\\

\textbf{$S$-step:} Use the representation $\V = \mathbf{X}^2\mathbf{S}\mathbf W^{2,\top}$ and test \eqref{eq:frictionFull} with $\mathbf{X}^2$ and $\mathbf{W}^2$. Then, the $S$-step equation reads
\begin{align}
    \dot{\mathbf{S}} = \mathbf{\widehat h}_2^{-2}\mathbf{S} \widehat{\SourceMatrix}_1^{\top}+\mathbf{\widehat h}_2^{-1}\mathbf{S} \widehat{\SourceMatrix}_2^{\top} + \mathbf{X}^{2,\top}\mathbf{u}_{m} \SourceVector^{\top}\mathbf{W}_2.
\end{align}
An implicit Euler time discretization gives with $\mathbf S^n = \mathbf{X}^{2,\top}\mathbf{X}^{1}\mathbf{S}^1\mathbf{W}^{1,\top}\mathbf{W}^2$
\begin{align}
    \mathbf{S}^{n+1} -\Delta t\mathbf{\widehat h}_2^{-2}\mathbf{S}^{n+1} \widehat{\SourceMatrix}_1^{\top}-\Delta t\mathbf{\widehat h}_2^{-1}\mathbf{S}^{n+1}\widehat{\SourceMatrix}_2^{\top}=  \mathbf{S}^{n}+ \Delta t\mathbf{X}^{2,\top}\mathbf{u}_{m}^{n+1,\top} \SourceVector^{\top}\mathbf{W}_2.
\end{align}
Hence, we need to solve an $r^2$ system of linear equations which requires $\mathcal{O}(r^6)$ operations. Moreover, computing flux matrices requires $\mathcal{O}((N_x+N)\cdot r^2)$ operations, hence $C_S^{f}\lesssim r^2\cdot (N+N_x) + r^6$.

\begin{remark}
Note that the implicit $L$-step exhibits high computational costs of $O(r^3\cdot N^3)$. Nevertheless, we are able to separate the spatial degrees of freedom from the number of moments, which are commonly much smaller than the number of spatial cells. Therefore, we expect a significant reduction of computational costs by DLRA, which, however, does not resemble its full potential. Note that DLRA will become substantially more efficient compared to the full-rank baseline for two-dimensional spatial domains since, in this case, the number of spatial cells drastically increases and thereby removes the computational bottleneck of matrix inversions performed in the $L$-step. In this case, the computational costs are distributed more evenly among the $K$, $L$, and $S$ steps.
\end{remark}

\section{Numerical experiments}
\label{sec:numerics}
As numerical examples of the model reduction techniques developed in the previous sections, we consider three test cases: (1) a water column (also called dam break) test case similar to the test case in \cite{Koellermeier2020} (2) a smooth wave similar to the test case in \cite{kowalski2018moment} and (3) the square root velocity profile test case that mimics a realistic velocity profile.
Additionally to POD-Galerkin and DLRA we also include computations of the HSWME with a reduced number of moments (rHSWME). In the following, $r$ denotes the rank in the case of DLRA and POD-Galerkin and the number of Moments for the rHSWME.
All numerical experiments can be reproduced with the openly available source code \cite{code}. The simulations have been performed on \texttt{11th Gen Intel(R) Core(TM) i7-11850H} CPUs (8 processing units).

\subsection{Water column or dam break test case}
\label{subsec:watercolumn}
In the first test case, we investigate a water column or dam break scenario defined by the settings in \cref{tab:setup_water_column}. Initially, the water is at rest and the height profile is defined as $h(x) = 0.3 + 0.35 \cdot \left(\text{tanh}(x)-\text{tanh}(x-0.2)\right)$, for $x\in[-1,1]$ leading to a water column within $[0,0.2]$ with slightly smoothed boundaries. Note, that we are using a smoothed water column, because of the oscillator behavior of a Lax-Friedrichs scheme with non-smooth initial conditions \cite{Breuss2004}. The friction parameters are $\nu = 1.0$ and $\lambda = 0.5$ which slow down velocities close to the bottom.
\begin{table}[htbp!]
    \centering
      \begin{tabular}{ll}
        \toprule
        friction coefficient & $\lambda=0.5$  \\
        slip length & $\nu=1.0$\\
        temporal domain & $t\in [ 0,0.2 ]$ \\
        spatial domain& $x\in [-1,1]$\\
        spatial resolution & $N_x=2000$ \\
        number of moments & $N=100$\\
        initial height & $h(x) = 0.3 + 0.35 \cdot \left(\text{tanh}(x)-\text{tanh}(x-0.2)\right)$ \\
        initial velocity & $u(0,x,\zeta)= 0 $ \\
        CFL number & $\text{CFL} = 0.25$\\
        spatial discretization & path-conservative FVM \cite{Pimentel2022}\\
        \bottomrule
      \end{tabular}
    \caption{Simulation setup for water column test case.}
    \label{tab:setup_water_column}
\end{table}

The numerical discretization is performed using $N_x=2000$ cells and a CFL number of $0.25$ for time stepping within $t\in [ 0,0.2 ]$. For the full-order HSWME model, we consider $N=100$ coefficients, leading to a large system of coupled PDEs. For the model reduction, we first consider a fixed number of $r=5$ basis function, which largely reduces the complexity of the reduced POD-Galerkin and DLRA models. 

\begin{figure}[htp!]
    \centering
    \setlength\figureheight{0.4\linewidth}%
    \setlength\figurewidth{0.6\linewidth}%
    \inputtikz{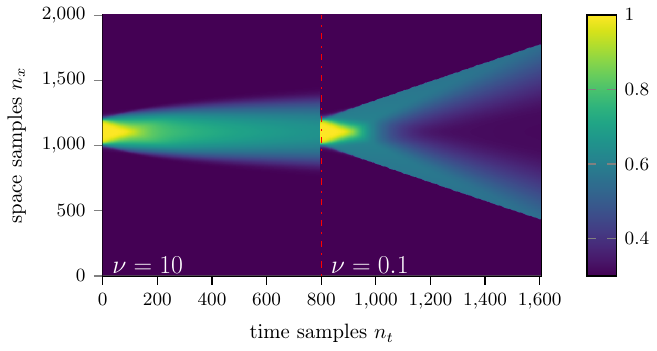}
    \caption{Sampled snapshot data of the POD, shown for the water height component.}
    \label{fig:water-column:snapshotdata}
\end{figure}

\begin{figure}[htp!]
    \centering
    \setlength\figureheight{0.4\linewidth}%
    \setlength\figurewidth{0.6\linewidth}%
    \inputtikz{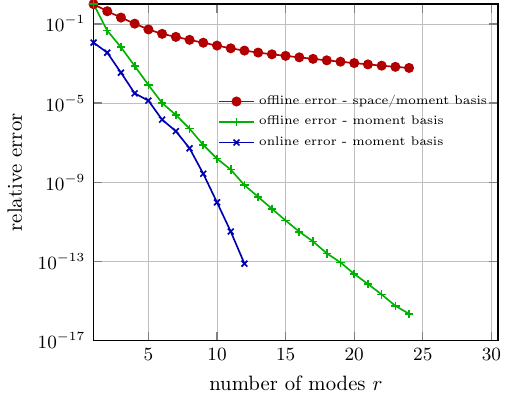}
    \caption{
        Relative approximation error $\sum_i\Vert{\vect{v}(x,t_i)-\tilde{\vect{v}}(x,t_i)}\Vert_{L_2}$
        of the snapshot data for space-moment $\tilde{\vect{v}}(x,t)=\sum_k \alpha_k(t) \vect{w}_k(x)$ 
        and moment only basis $\tilde{\vect{v}}(x,t) =\sum_{k=1}^r \hat{\alpha}_k(x,t) \vect{w}_k$. The online error is given as $\Vert{\vect{u}(x,t_\text{end})-\tilde{\vect{u}}(x,t_\text{end})}\Vert_{L_2}$, where $\vect{u}=(h, h u_m)$ 
    }
    \label{fig:water-column:error-decay-pod}
\end{figure}

The trial and test space for the POD-Galerkin approach is spanned by the POD-basis. It is setup offline from the concatenation of the snapshots for two different trajectories sampled with slip length $\nu\in\{0.1,10\}$. The snapshots are shown for the water height in \cref{fig:water-column:snapshotdata}. From the snapshot matrix $\matr{V}^\mathrm{POD}\in \mathbb{R}^{(2N_tN_x)\times N}$ we can compute the POD basis with help of the SVD \cref{eq:svd-POD} and evaluate the projected right-hand side efficiently for the testing slip length $\nu=1$.
We highlight that our approach does not try to decouple space and time and we therefore avoid slowly decaying approximation errors of our dyadic decomposition, a known problem for transport-dominated fluid systems. The difference between the classical MOR approach, which suffers from this slow decay and the here presented approach is shown in \cref{fig:water-column:error-decay-pod}. As our approach only tries to reduce the component space, it avoids slowly decaying approximation errors, thus yielding a rapid decay of the POD approximation errors. Note, that this procedure is nothing new but has been extensively used in hierarchical MOR reduction (see for example \cite{SmentanaOhlberger2017}).

\begin{figure}[htp!]
    \centering
    \begin{subfigure}{0.45\linewidth}
        \centering
        \setlength\figureheight{0.9\linewidth}%
        \setlength\figurewidth{1\linewidth}%
        \caption{}
        \vspace{-0.8cm}
        \inputtikz{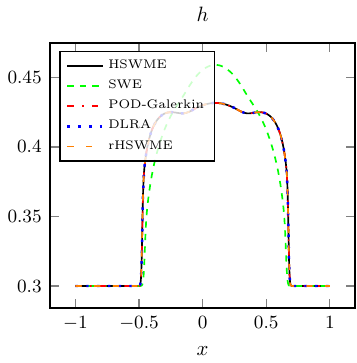}%
    \end{subfigure}%
    \begin{subfigure}{0.45\linewidth}
        \centering
        \setlength\figureheight{0.9\linewidth}%
        \setlength\figurewidth{1\linewidth}%
        \caption{}
        \vspace{-0.8cm}
        \inputtikz{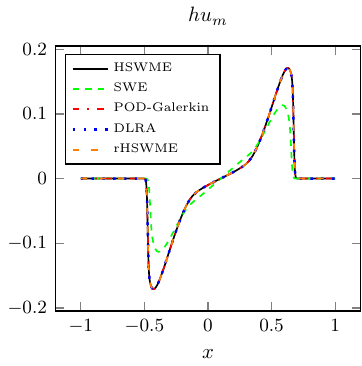}%
    \end{subfigure}
    \caption{Water column test case comparison of macroscopic quantities water height $h$ (a) and momentum $h u_m$ (b) for full-order HSWME, SWE, POD-Galerkin, and DLRA using rank $r=5$. Both reduced models DLRA and POD-Galerkin yield indistinguishable results from the full-order model while the simple SWE model shows insufficient accuracy.}
    \label{fig:water-column:state-profiles}
\end{figure}

The numerical results shown in \cref{fig:water-column:state-profiles} compare the numerical results of (1) the full-order HSWME model, (2) the lowest-order SWE model, (3) the new macro-micro decomposition POD-Galerkin, (4) the new macro-micro decomposition DLRA and (5) the HSWME with a reduced number of moments. The full-order HSWME model clearly depicts two waves symmetrically moving left and right, respectively. The simple SWE model uses a constant velocity profile with zero micro structure and thus cannot capture the complex dynamics induced by the bottom friction, which slows down the velocity profile at the bottom. The new macro-micro decomposition POD-Galerkin and DLRA methods, yield optically indistinguishable solutions from the full-order HSWME model, even though only $r=5$ basis functions are chosen. This emphasizes the good approximation quality of the reduced models despite a small number of basis functions. However, we note that for this specific test case the moment model with only 5 moments already does a good job as well. This is because the dependencies on the higher moments vanishes quickly and a small number of moments already gives an excellent approximation quality.

\begin{figure}
    \centering
    \begin{subfigure}{0.9\linewidth}
        \centering
        \setlength\figureheight{0.45\linewidth}%
        \setlength\figurewidth{0.7\linewidth}%
        \caption{}
        \inputtikz{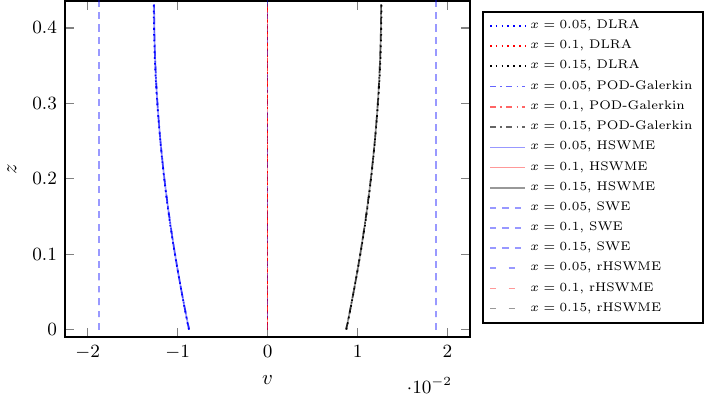}%
        \label{fig:water-column:velocity-profiles_close}
    \end{subfigure}
    \begin{subfigure}{0.9\linewidth}
        \centering
        \setlength\figureheight{0.45\linewidth}%
        \setlength\figurewidth{0.7\linewidth}%
        \caption{}
        \inputtikz{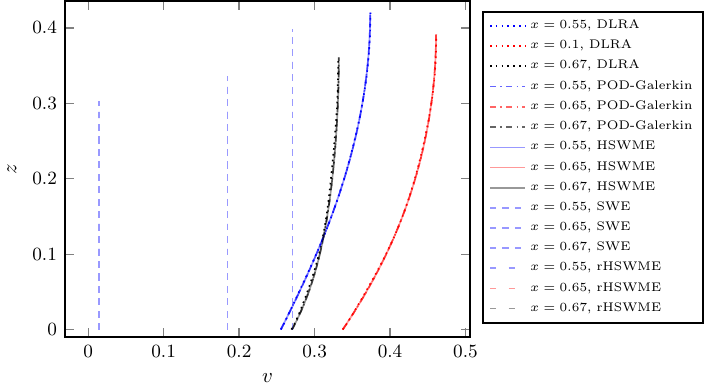}%
        \label{fig:water-column:velocity-profiles_far}
    \end{subfigure}
    \caption{Water column test case comparison of velocity profiles close to the domain center (a) and close to the shock wave (b) for full-order HSWME, SWE, POD-Galerkin, and DLRA using rank $r=5$. Both reduced models DLRA and POD-Galerkin also yield indistinguishable results from the full-order model while the simple SWE model shows insufficient accuracy and wrong propagation speeds.}
    \label{fig:water-column:velocity-profiles}
\end{figure}

In \cref{fig:water-column:velocity-profiles} the water velocity profiles are plotted each time for different positions: $x = 0.05$, 0.0, 0.15 close to the center of the domain in \cref{fig:water-column:velocity-profiles_close} and $x = 0.65, 0.67$ close to the shock wave in \cref{fig:water-column:velocity-profiles_far}. In \cref{fig:water-column:velocity-profiles_close}, we clearly see that also the velocity profiles of the full-order model and the reduced POD-Galerkin and DLRA models, as well as the reduced moment model agree almost perfectly at all three points. The SWE model on the other hand shows an overestimation of the average velocity at both positions. This is due to wrong propagation speeds of the SWE model \cite{Koellermeier2020} and clearly shows why this simple model is not useful in simulations of such complex cases. In \cref{fig:water-column:velocity-profiles_far}, the profiles close to the shock wave are plotted and the reduced models again agree almost perfectly with the full-order model. Again the SWE model predicts wrong average profiles, in this case much smaller than the full model and the reduced models.

\begin{figure}
    \centering
    \setlength\figureheight{0.5\linewidth}%
    \setlength\figurewidth{0.6\linewidth}%
    \inputtikz{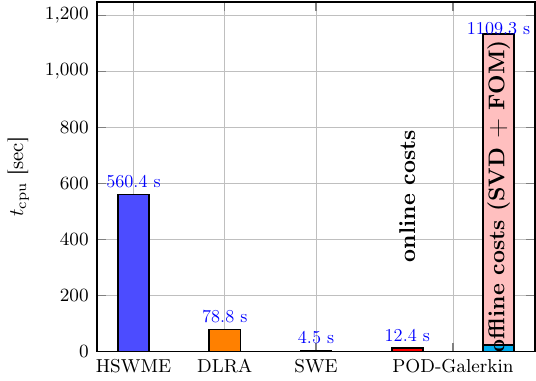}
    \caption{Water column test case runtime comparison between full-order HSWME, DLRA, SWE, POD-Galerkin (including offline precomputations and online phase separately) show a large speedup of the reduced models with respect to the full model. The simple SWE model is obviously fastest but does not achieve sufficient accuracy. The truncation ranks of DLRA $r=4$ and POD-Galerkin $r=3$ are tuned such that the relative $L_2$ errors of $\boldsymbol{u}$ at the final time are approximately the same ($\approx 0.3\%$).}
    \label{fig:water-column:performance_r}
\end{figure}

\begin{figure}
    \centering
    \setlength\figureheight{0.5\linewidth}%
    \setlength\figurewidth{0.75\linewidth}%
    \inputtikz{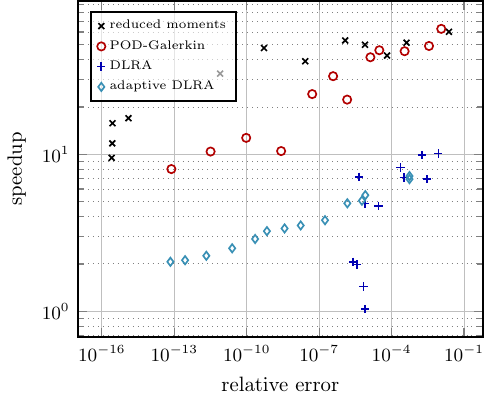}%
    \caption{Water column test case runtime comparison speedup and error comparison between reduced moments HSWME, DLRA, POD-Galerkin (without offline precomputation phase) in comparison to full HSWME model \refereeOne{(100 moments)}.
     \refereeOne{For the comparison, we vary the degrees of freedom in the system by adjusting the rank (POD-Galerkin, DLRA) or the number of moments (reduced moments) to $r=1,2,3,4,5,7,9,11,15,20,25,30$. For the rank adaptive DLRA, we select 14 different tolerance values on a logarithmic scale from $10^{-10}$ to $10^0$.} 
     The relative error denotes the $L_2$ error of $u=(h,hu_m)$ at the final time $t=0.2$.}
    \label{fig:water-column:performance}
\end{figure}

For the settings from above, \cref{fig:water-column:performance_r} shows the runtime comparison between the full-order HSWME, DLRA, SWE, and POD-Galerkin, where the POD-Galerkin runtime is divided into the offline precomputations and the online phase.
The truncation ranks of DLRA $r=4$ and POD-Galerkin $r=3$ are tuned such that the relative $L_2$ errors of $\boldsymbol{u}$ at the final time are approximately the same ($\approx 0.3\%$, compare \cref{fig:water-column:performance}). The DLRA method already reduces the runtime significantly. The SWE has the fastest runtime, but does not result in sufficient accuracy as seen in the previous figures. The POD-Galerkin method is equally as fast as the SWE during the online phase, but requires a relatively costly offline precomputation phase.
It is important to highlight that the efficacy of a global reduced order model becomes particularly apparent in the realm of multi-query simulations. In this context, the upfront costs associated with the offline phase of the POD-Galerkin algorithm are offset by an exceptionally efficient online phase, which can be replicated across a considerable number of queries.

However, there might be other MOR applications in which a global ROM can not be set up, because of the tremendous amount of data that would be required in the offline stage. \refereeOne{For example, in the realistic scenario of two-dimensional HSWME or higher-dimensional systems involving multiple parameters formulated in a tensor-valued ODE (discussed in \cite{OthmarLubich2010}).}
This is where DLRA has a clear advantage, as it does not rely on an offline phase. \refereeOne{When the offline phase becomes too expensive, a combination of both methods might be beneficial. In this approach, DLRA is used to establish the initial basis, which is then employed in the online phase of POD.}

Increasing the rank of the reduced models leads to increasing accuracy at the expense of more runtime, as can be seen in \cref{fig:water-column:performance}, where we compare the online computation time of POD-Galerkin with DLRA and the reduced moment HSWME model. For the rank adaptive version of DLRA we sample the solution for 14 different $\vartheta$ values evenly on a log scale in the interval $10^{-10}\le\vartheta\le10^0$.  From the comparison we observe, that rHSWME is the fastest, which is not surprising since the higher moments vanish quickly, in a regime of high friction. POD-Galerkin seems to achieve similar errors for small number of modes, however, it is slower then the rHSWME with increasing accuracy. Both DLRA and its adaptive version are one order of magnitude slower than POD-Galerkin. 
Note that since we apply MOR only in the microscopic higher moments ($v$-components) the reduced DLRA model converges towards the SWE with $r\to0$. Therefore, we obtain small relative errors even with a small number of modes.

As a summary of this first test case, we see very good accuracy of the reduced models with small rank and significant speedup that can even be amplified by further reducing the rank. Unfortunately, these speedups do not pay off directly when comparing it to the reduced moments HSWME, as the higher moments vanish quickly, due to the choice of the test case. In the next section, we will see, that this is not generally the case.

\subsection{Smooth wave}\label{sec:smoothWave}
This second test case  closely follows the general simulation test cases in \cite{kowalski2018moment} and its setup is given in \cref{tab:setup_smooth_wave}. It describes a smooth wave given by the initial height function $h(x) = 1 + \text{exp}(3 \cos(\pi (x + 0.5)))/\text{exp}(4)$ travelling through a periodic domain $x\in [-1,1]$. 
The initial velocity profile is chosen as $u(0,x,\zeta)= 0.25\cdot (1 - \phi_1(\zeta) + \phi_N(\zeta))$, leading to $u_m = 0.25$, $\alpha_1 = -0.25$, and $\alpha_N=0.25$. This velocity profile happens to be outside of the hyperbolicity region of the standard SWME model in the test case in \cite{Koellermeier2020}. This means that applying the SWME model can lead to instability issues. This is resolved by applying the HSWME model and demonstrates the utility of this guaranteed hyperbolic model necessary to result in a well-posed model.
Given the substantial dependence of the initial profile on the last moment $\alpha_N$, it is anticipated that the rHSWME will be surpassed by DLRA and POD-Galerkin.

\begin{table}[htbp!]
    \centering
      \begin{tabular}{ll}
        \toprule
        friction coefficient & $\lambda=0.001$  \\ 
        slip length & $\nu=100$\\        
        temporal domain & $t\in [ 0,0.2 ]$ \\
        spatial domain& periodic $x\in [-1,1]$\\
        spatial resolution & $N_x=2000$ \\
        number of moments & $N=100$ \\
        initial height & $h(x) = 1 + \text{exp}(3 \cos(\pi (x + 0.5)))/\text{exp}(4)$ \\
        initial velocity & $u(0,x,\zeta)= 0.25\cdot (1 - \phi_1(\zeta) + \phi_N(\zeta))$ \\
        CFL number & $\text{CFL} = 0.2$\\
        spatial discretization & path-conservative FVM \cite{Pimentel2022}\\
        \bottomrule
      \end{tabular}
    \caption{Simulation setup for smooth wave test case.}
    \label{tab:setup_smooth_wave}
\end{table}

The numerical discretization uses the same parameters as before. The test is carried out using $N_x=2000$ cells and a CFL number of $0.25$ for the time interval $t\in [ 0,0.2 ]$. The full-order HSWME model uses $N=100$ coefficients and the reduced models use a fixed number of $r=4$ basis function.

Similar to the water column test case, we compute the POD-basis from the snapshots of two trajectories simulated with the full HSWME at $\nu\in\{10,1000\}$. The friction parameters to compare all methods are now chosen as $\nu = 10$ and $\lambda = 0.001$, leading to large values of the bottom friction term.

\begin{figure}
    \centering
    \begin{subfigure}{0.45\linewidth}
        \centering
        \setlength\figureheight{0.9\linewidth}%
        \setlength\figurewidth{1\linewidth}%
        \caption{}
        \vspace{-0.8cm}
        \inputtikz{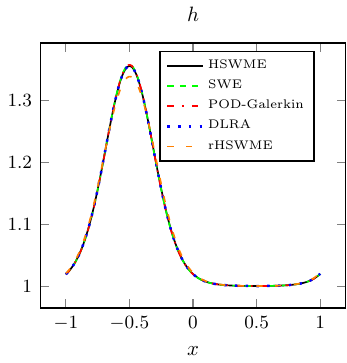}%
    \end{subfigure}%
    \begin{subfigure}{0.45\linewidth}
        \centering
        \setlength\figureheight{0.9\linewidth}%
        \setlength\figurewidth{1\linewidth}%
        \caption{}
        \vspace{-0.8cm}
        \inputtikz{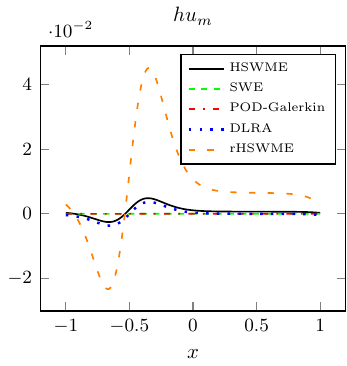}%
        \label{subfig:KT:hu-profile}
    \end{subfigure}
    \caption{Smooth wave test case comparison of macroscopic quantities water height $h$ 
    (a) and momentum $h u_m$ 
    (b) for full and reduced-moments HSWME, SWE, POD-Galerkin, and DLRA using rank $r=4$. Again, the reduced models DLRA and POD-Galerkin yield very good approximations of the full-order model while the simple SWE model fails.}
    \label{fig:KT:state-profiles}
\end{figure}

We first show the numerical results for the macroscopic water height $h$ and momentum $h u_m$ in \cref{fig:KT:state-profiles} for (1) the full-order HSWME model, (2) the lowest-order SWE model, (3) the reduced HSWME model, (4) the new macro-micro decomposition POD-Galerkin, and (5) the new macro-micro decomposition DLRA. As for the first test case, the simple SWE model fails at capturing the complex dynamics of the test case, which is especially apparent in the momentum $h u_m$ shown in \cref{subfig:KT:hu-profile}. Similarly, POD Galerkin and the rHSWME model struggle to capture an accurate description of the momentum. Here DLRA yields a very good match even with a small number of degrees of freedom.

\begin{figure}
    \centering
    \begin{subfigure}{0.9\linewidth}
        \centering
        \setlength\figureheight{0.45\linewidth}%
        \setlength\figurewidth{0.7\linewidth}%
        \caption{}
        \inputtikz{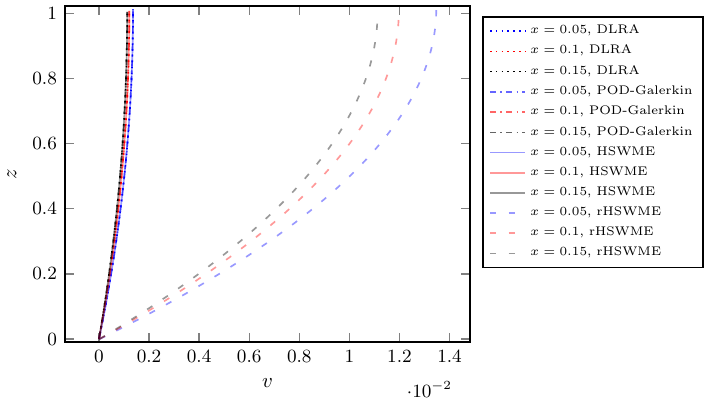}%
        \label{fig:KT:velocity-profiles_close}
    \end{subfigure}
    \begin{subfigure}{0.9\linewidth}
        \centering
        \setlength\figureheight{0.45\linewidth}%
        \setlength\figurewidth{0.7\linewidth}%
        \caption{}
        \inputtikz{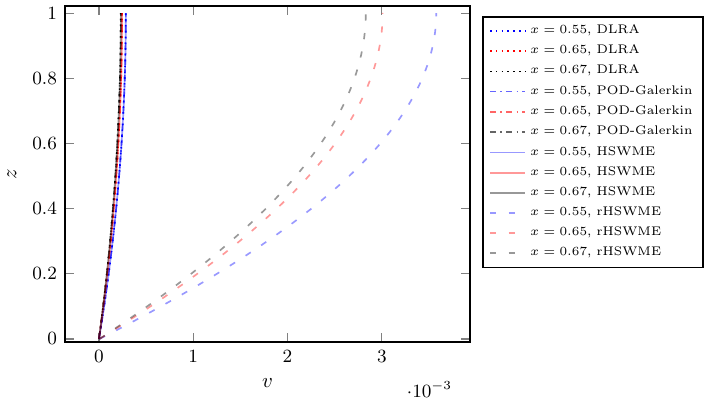}%
        \label{fig:KT:velocity-profiles_far}
    \end{subfigure}
    \caption{Smooth wave test case comparison of velocity profiles close to the domain center 
    (a) and close to the shock wave 
    (b) for full-order HSWME, POD-Galerkin, and DLRA using rank $r=5$ (SWE omitted for conciseness). Both reduced models DLRA and POD-Galerkin yield good results in comparison with the full-order model.}
    \label{fig:KT:velocity-profiles}
\end{figure}

In \cref{fig:KT:velocity-profiles} the water velocity profiles for the smooth wave test case are plotted each time for different positions: $x = 0.05, 0.0, 0.15$ close to the center of the domain in \cref{fig:KT:velocity-profiles_close} and $x = 0.65, 0.67$ further away from the center in \cref{fig:KT:velocity-profiles_far}. In \cref{fig:KT:velocity-profiles_close}, the velocity profiles of the full-order model agree well with the reduced POD-Galerkin and DLRA models at all three points, while there are some small differences in the maximum velocity value due to the complexity of the test case and the significant model reduction. In contrast the rHSWME completely overshoots in the maximum velocity as it misses the information of the highest moment. Similarly for the SWE model, that is not shown for conciseness. 
The velocity profiles at positions further away from the center in  \cref{fig:KT:velocity-profiles_far} yield a similar result. 

\begin{figure}
    \centering
    \setlength\figureheight{0.25\linewidth}%
    \setlength\figurewidth{0.8\linewidth}%
    \inputtikz{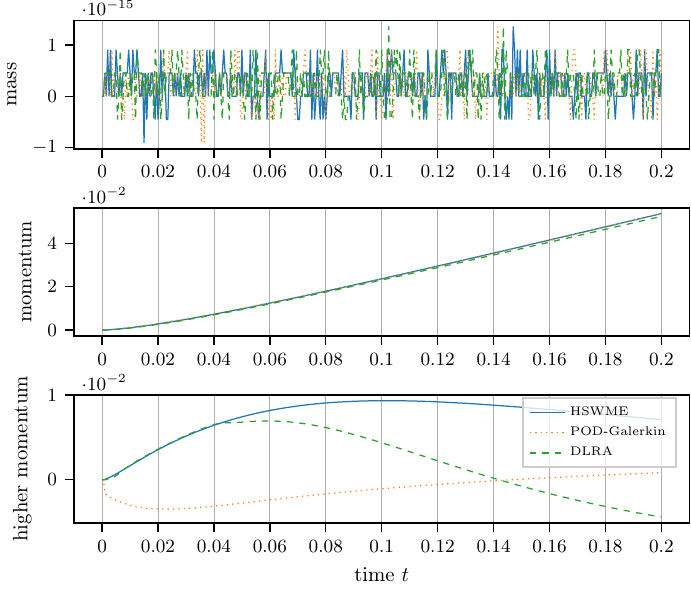}%
    \caption{Smooth wave test case time evolution of the total mass ($\int_{-1}^1 h\, \mathrm{d} x$), momentum ($\int_{-1}^1hu\,\mathrm{d} x$) and higher momentum ($\int_{-1}^1h\alpha_4\,\mathrm{d} x$) relative to their respective initial values plotted for the full HSWME model and both reduced POD-Galerkin and DLRA models. The mass is conserved and the momentum accurately follows the full model reference solution.}
    \label{fig:water-column:conserv}
\end{figure}

For this smooth wave test case, we also want to emphasize one main property of our newly developed macro-micro decomposition reduced models, which is guaranteed mass conservation by construction. In \cref{fig:water-column:conserv}, the time evolutions of the total mass $\int_{-1}^1 h\, \mathrm{d} x$, the total momentum $\int_{-1}^1 h u_m\, \mathrm{d} x$ and the total fourth higher moment $\int_{-1}^1 h \alpha_4\, \mathrm{d} x$ are plotted in terms of relative deviation from the initial values. It is clearly seen that the total mass is constant as changes are within machine precision. The total momentum and total higher momentum are not conserved in agreement with the underlying PDE model, which includes friction terms for the corresponding equations leading to a loss of momentum due to bottom friction, for example. The mass conservation is achieved due to the macro-micro decomposition formulation which includes explicitly solving for the water height while applying the model reduction only to the remaining microscopic velocity profile coefficients. While the total momentum is not constant, its evolution is in very good agreement with the full HSWME model for both reduced models since also the equation for $h u_m$ is apart from the model reduction process. Some deviations are seen for the total fourth momentum, obviously originating from the model reduction process for that variable.

\begin{figure}
    \centering
    \setlength\figureheight{0.25\linewidth}%
    \setlength\figurewidth{0.8\linewidth}%
    \inputtikz{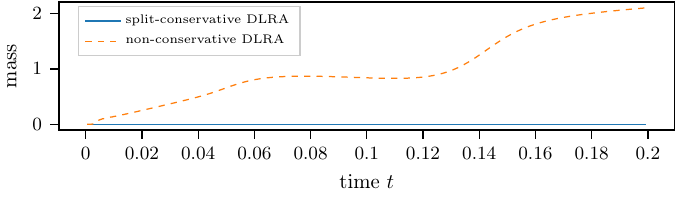}%
    \caption{Smooth wave test case time evolution of the total mass ($\int_{-1}^1 h\, \mathrm{d} x$) relative to its respective initial value plotted for the novel macro-micro decomposition conservative DLRA and a naive (non-conservative) DLRA. The mass is only conserved for the macro-micro decomposition conservative DLRA.}
    \label{fig:water-column:non-conserv}
\end{figure}

With \cref{fig:water-column:non-conserv} we want to emphasize that a naive application of model reduction techniques does not lead to conservation of mass. This is done by comparison of our new macro-micro decomposition conservative DLRA, where the evolution of the macroscopic water height $h$ and momentum $h u_m$ is decoupled from the microscopic reduced coefficient system, with a naive (non-conservative) DLRA, where the complete system including water height $h$, momentum $h u_m$ and coefficients $h \alpha_i$ is reduced as a whole. It is seen that only the macro-micro decomposition conservative DLRA method achieves mass conservation.

\begin{figure}[htp!]
    \centering
    \setlength\figureheight{0.25\linewidth}%
    \setlength\figurewidth{0.8\linewidth}%
    \inputtikz{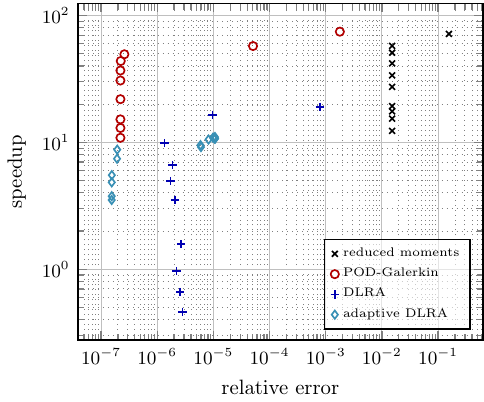}%
    \caption{Smooth wave test case runtime comparison speedup and error comparison between rHSWME, SWE, DLRA, POD-Galerkin (without offline precomputation phase) in comparison to full HSWME model \refereeOne{(100 moments)}.
    For the comparison, we vary the degrees of freedom in the system by adjusting the rank (POD-Galerkin, DLRA) or the number of moments (reduced moments) to $r=1,2,3,5,7,9,15,20,25,30$. For the rank adaptive DLRA, we select 14 different tolerance values on a logarithmic scale from $10^{-10}$ to $10^0$. 
    The relative error denotes the $L_2$ error of $u=(h,hu_m)$ at the final time $t=0.2$. }
    \label{fig:smooth:speed-vs-error}
\end{figure}

Next, we compare the runtime vs. speedup of the presented methods. For this study, the rank of DLRA and POD-Galerkin is gradually increased starting from $r=1$ until $r=30$ (not all points in between are included).  As before the rank adaptive version of DLRA is sampled for 14 different $\vartheta$ values evenly on a log scale in the interval $10^{-10}\le\vartheta\le10^0$.
Comparing both MOR methods and the rHSWME in \cref{fig:smooth:speed-vs-error}, we observe that rHSWME stagnates at about 2\% error. This can be explained as after the linear moment the higher moments only increase the complexity of the calculations, but can not capture any sensible dynamics before the last moment is not included. DLRA and POD-Galerkin however do not a-priori impose a basis and therefore can include the information necessary to represent the dynamics of the highest moments. As expected, we observe that DLRA exhibits an increased runtime in comparison to the online phase of POD-Galerkin. In general, we see that the adaptive DLRA method can achieve results with higher accuracy, which is compromised by additional effort for adding and removing basis functions to the decomposition. However, DLRA does not require a computationally expensive and memory-intensive offline phase. Furthermore, since DLRA utilizes time-dependent basis functions, it exhibits two advantages over POD-Galerkin. First, since basis information can be added and removed in time, the approximation space at a given rank is richer than for POD-Galerkin. This is for example seen when comparing both methods for a fixed rank $r=4$ in \cref{fig:KT:state-profiles}. DLRA achieves a better approximation quality with fewer degrees of freedom.
Second, when the dynamics in the online phase is not captured by the POD-ansatz space, POD-Galerkin requires an expensive retraining while DLRA adapts automatically to such situations. A main advantage of POD-Galerkin, besides the reduced runtime is, that it does not require the derivation of new evolution equations. Moreover, in the case of certain non-linearities, efficient evolution equations for DLRA might not be available.

\subsection{Square root profile}\label{sec:SquareRoot}

In this section, we present a realistic test case that implements a water column with an initial square root velocity profile $u(0,x,\zeta)= \sqrt{\zeta}$ in the viscous case. Note, that in comparison to \cref{subsec:watercolumn} we use an increased slip length and decreased friction coefficient that highlights a stronger coupling of the higher moments to $(h,hu_m)$. The range of parameters is similar to the test cases used in \cite{Huang2022}, where the relaxation of the model towards different equilibria is investigated.
The precise test setup is detailed in \cref{tab:setup_sqrt}.

\begin{table}[htp!]
    \centering
      \begin{tabular}{ll}
        \toprule
        friction coefficient & $\lambda=0.01$  \\
        slip length & $\nu=10$\\
        temporal domain & $t\in [ 0,0.05 ]$ \\
        spatial domain& periodic $x\in [-0.15,0.3]$\\
        spatial resolution & $N_x=2000$ \\
        \refereeOne{number of moments} & \refereeOne{$N=100$} \\
        initial height & $h(x)=\frac{0.7}{2}(-\tanh(50(x - 0.2)) + \tanh(50(x))) + 0.3$ \\
        initial velocity & $u(0,x,\zeta)= \sqrt{\zeta}$ \\
        CFL number & $\text{CFL} = 0.1$\\
        spatial discretization & path-conservative FVM \cite{Pimentel2022}\\
        \bottomrule
      \end{tabular}
    \caption{Simulation setup for square root profile test case, compare \cite{Pimentel2022}.}
    \label{tab:setup_sqrt}
\end{table}

As before we sample the solution at two different slip lengths $\nu=1$ and $\nu=100$ to set up the POD basis in an offline phase. For each parameter, we sample 800 snapshots, resulting in a total of 1600 snapshots to set up the POD basis as detailed in \cref{subsec:offlinephase}.

\begin{figure}
    \centering
    \begin{subfigure}{0.45\linewidth}
        \centering
        \setlength\figureheight{0.8\linewidth}%
        \setlength\figurewidth{1\linewidth}%
        \caption{}
    \inputtikz{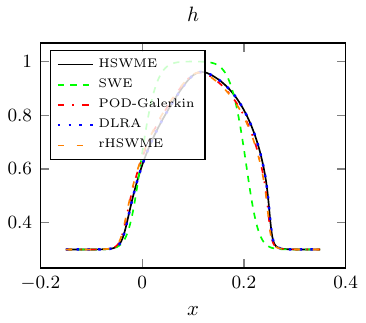}%
        \label{subfig:sqrt:h}
    \end{subfigure}
    \begin{subfigure}{0.45\linewidth}
        \centering
        \setlength\figureheight{0.8\linewidth}%
        \setlength\figurewidth{1\linewidth}%
        \caption{}
    \inputtikz{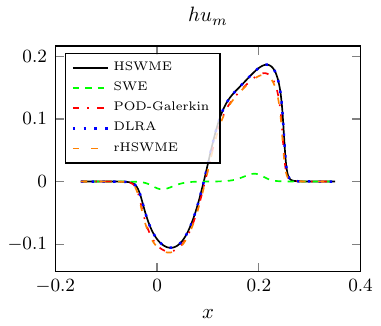}%
        \label{subfig:sqrt:hu}
    \end{subfigure}
    \caption{Square root test case comparison of macroscopic variables $h$ (a) and $h u_m$ (b) for full-order HSWME, SWE, POD-Galerkin, and DLRA using rank $r=4$ and a low-order rHSWME model with only 4 moments. DLRA shows excellent accuracy and POD-Galerkin yield good results in comparison with the full-order model.}
    \label{fig:sqrt-profiles}
\end{figure}

The simulation results are shown in \cref{subfig:sqrt:h} and \cref{subfig:sqrt:hu}. We see a very good match of the DLRA method with the full HSWME model, while POD-Galerkin and the rHWSME model (with less coefficients $r=4$) shows a small deviation from the full model solution. This square root profile example shows that DLRA and its adaptive version are superior for a small number of degrees of freedom. 
This is especially relevant for realistic scenarios in higher space dimension because here one does not have the luxury of choosing a high number of degrees of freedom due to the memory bottleneck.

\begin{figure}
    \centering
    \begin{subfigure}{0.9\linewidth}
        \centering
        \caption{}
        \setlength\figureheight{0.8\linewidth}%
        \setlength\figurewidth{1\linewidth}%
        \inputtikz{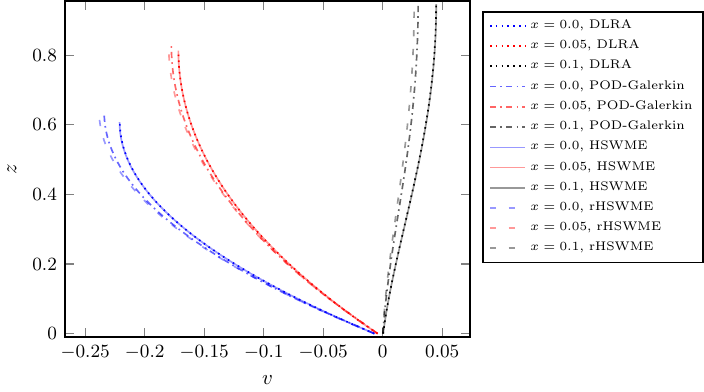}%
        \label{subfig:sqrt:velocity-profiles_close}
    \end{subfigure}
    \begin{subfigure}{0.9\linewidth}
        \centering
        \setlength\figureheight{0.8\linewidth}%
        \setlength\figurewidth{1\linewidth}%
        \caption{}
        \inputtikz{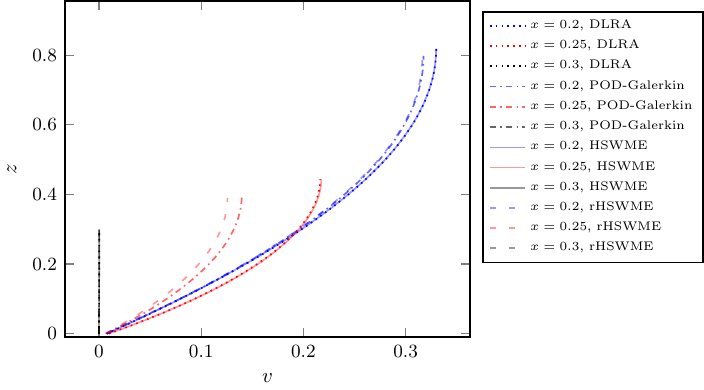}%
        \label{subfig:sqrt:velocity-profiles_far}
    \end{subfigure}
    \caption{Square root test case comparison of velocity profiles close to the domain center (a) and close to the shock wave (b) for full-order HSWME, POD-Galerkin, DLRA using rank $r=4$ and a low-order rHSWME model with only $r=4$ modes (SWE omitted for conciseness). Again DLRA shows excellent accuracy and POD-Galerkin yield good results in comparison with the full-order model.}
    \label{fig:sqrt-v-profiles}
\end{figure}

\Cref{fig:sqrt-v-profiles} also visualizes the difference in the velocity profiles. The same observation as for the macroscopic variables $h, u_m$ can be made: DLRA succeeds at reproducing the full model solution of HSWME, while POD-Galerkin and the smaller rHWSME model lead to an error in the velocity profiles. This is true both for a point in the center of the simulation domain in \cref{subfig:sqrt:velocity-profiles_close} as well as for a point close to the shock wave in \cref{subfig:sqrt:velocity-profiles_far}, where the square root velocity profile is still clearly visible.

\begin{figure}[htp!]
    \centering
    \setlength\figureheight{0.25\linewidth}%
    \setlength\figurewidth{0.8\linewidth}%
    \inputtikz{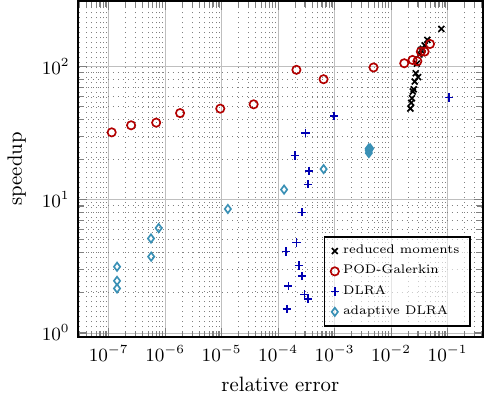}%
    \caption{Square root test case runtime comparison speedup and error comparison between DLRA, POD-Galerkin (without offline precomputation phase) and rHSWME in comparison to full HSWME model \refereeOne{(100 moments)}. 
      \refereeOne{For the comparison, we vary the degrees of freedom in the system by adjusting the rank (POD-Galerkin, DLRA) or the number of moments (reduced moments) to $r=1,2,3,4,5,6,7,8,9,10,11,15,20,25,30$. For the rank adaptive DLRA, we select 14 different tolerance values on a logarithmic scale from $10^{-10}$ to $10^0$.
      }The relative error denotes the $L_2$ error of $(h,hu_m)$ at the final time $t=0.05$.
}
    \label{fig:sqrt:speed-vs-error}
\end{figure}

\Cref{fig:sqrt:speed-vs-error} shows the speedup obtained depending on the resulting error simulated with the respective method. Each data point represents one choice for the number of the reduced model coefficients $r$ or $10^{-10}\le\vartheta\le 10^0$. Note that offline computation costs are neglected here. As expected and shown for the previous test cases, POD-Galerkin has a very fast online computation phase leading to large speedups. The rank-adaptive DLRA can obtain very good accuracy as well even without any offline computation, but at the sacrifice of a more costly online computation. Simply reducing the number of moments in the full model (called rHSWME above) leads to significant errors and is therefore not recommended.

The square root test demonstrates for a more physically relevant velocity profile the power of the rank-adaptive model order reduction via adaptive DLRA.  

\subsection{Summary of numerical results}
\refereeThree{To summarize and compare the gained insights from numerical test cases, we note that, as expected, reduced moment models (rHSWME) perform well in settings with high friction and smooth initial velocity profiles. This is especially seen in Section~\ref{subsec:watercolumn}, where rHSWME yields the best performance, achieving a speedup of over 50 compared to the full model to achieve a relative error smaller than $10^{-5}$. While POD Galerkin demonstrates a speedup of around 20 at the same error level, the additional computational costs required by DLRA lead to a significantly lower speedup of around 5. In situations where the initial velocity profile is irregular or friction is small, POD Galerkin and DLRA exhibit a clear advantage over rHSWME. This is seen in Sections~\ref{sec:smoothWave} and \ref{sec:SquareRoot}, where an error level of $10^{-5}$ is surpassed at a speedup of around 50 for POD Galerkin and 16 for DLRA in Section~\ref{sec:smoothWave}, and 50 for POD Galerkin and 8 for DLRA in Section~\ref{sec:SquareRoot}. Overall, POD Galerkin achieves a higher speedup at given error levels compared to DLRA when factoring out offline costs. This is expected, as in this case, POD Galerkin uses a fixed basis, whereas DLRA additionally requires solving basis update equations. However, POD-Galerkin requires additional hyperparameters, while DLRA allows errors to be controlled during computation when using rank-adaptive integrators with a single hyperparameter. Comparing fixed-rank and rank-adaptive DLRA methods, it is observed that an adaptive rank will significantly improve the method's performance. The two factors that benefit the rank-adaptive method are the exact initial conditions of the $S$-steps and the ability to pick varying ranks over time.}

\section{Conclusion}
\label{sec:conclusion}
In this work we proposed mass conservative model order reduction methods for the hyperbolic shallow water moment equations that yield fast and accurate solutions. Mass conservation is achieved by decomposing the macroscopic water height and momentum equations from the microscopic higher-order moments and applying model order reduction solely to the microscopic higher-order moment equations. Additionally, the decomposition allows to recover the naive SWE for vanishing truncation rank of the reduced HSWME, resulting in a very accurate reduced model even for a small number of modes. We use two model order reduction methods, namely POD-Galerkin and dynamical low-rank approximations to accelerate the computation and reduce memory footprint. The methods can produce speedups of up to 100 compared to the full HSWME, while the introduced approximation errors are negligible. However, it must be noted, that the HSWME with a reduced number of moments already shows a good performance, if the higher moments vanish quickly. Concluding that DLRA and POD-Galerkin are most efficient if the shearing in the vertical direction is strong and therefore a slow decay in the additional moments is expected. Employing rank-adaptive methods led to further runtime improvements while achieving high accuracy with a smaller number of variables.

The work in this paper opens up possibilities for future work on model reduction for shallow water moment models. An interesting extension would be to use DLRA to generate basis functions for POD.

\section*{Author Contribution Statement (CRediT)}
All authors contributed equally to this publication. The original draft and all subsequent edits and reviews were done equally by all authors.
The authors contributions differ in the following points:

\vspace{5pt}
{\small
\noindent
\begin{tabular}{@{}lp{10.8cm}}
\textbf{Julian Koellermeier:} & initial idea,  methodology, implementation of HSWME, setup of test cases\\
\textbf{Philipp Krah:} & implementation of POD-Galerkin, simulation/setup of numerical tests,  visualization\\
\textbf{Jonas Kusch:} & initial idea, methodology, implementation of HSWME, implementation of DLRA, simulation/setup of numerical tests, visualization\\
\end{tabular}
}

\section*{Conflict of Interest}
The authors declare that they have no conflict of interest.

\section*{Code and Data Availability}
All scripts to reproduce the results are available at \cite{code}.

\section*{Acknowledgement}
The authors would like to acknowledge the financial support of the CogniGron research center and the Ubbo Emmius Funds (University of Groningen). Jonas Kusch was funded by the Deutsche Forschungsgemeinschaft (DFG, German Research Foundation) – 491976834.
The authors were granted access to the HPC resources of IDRIS under the allocation No. AD012A01664R1 attributed by Grand Équipement National de Calcul Intensif (GENCI).
Centre de Calcul Intensif d’Aix-Marseille is acknowledged for granting access to its high performance computing resources.
Furthermore, we would like to thank Yannick Schubert for detailed reading and commenting on the work.

\bibliographystyle{abbrv}
\bibliography{references}

\appendix
\section{Numerical discretization of the decomposed transport part}
\label{app1}
This appendix briefly describes a numerical scheme for the transport part and then continues with the space-time discretization for the decomposed set of variables \cref{eq-def:statematrixV}.

We first cover a standard scheme.
For the spatial discretization, we consider an equidistant mesh $x_j = j\Delta x $ with lattice spacing $\Delta x$ and local state-variables, so that the full state vector at time $t$ is given by $\Q(t) = [ \vect{q}_i(x_j,t) ]_{i j} \in \mathbb{R}^{N_x \times (N+2)}$ for $N_x$ spatial discretization points.
Furthermore, we denote the rows of this matrix as $\Q_j(t)\in\mathbb{R}^{N+2}$, which represents the corresponding solution vector evaluated at the spatial cell $j$.

The spatial discretization of \cref{e:transport} for each cell $j$ leads to a semi-discrete ODE for the unknowns $\Q_j(t)$. This step is performed via the following standard first-order path-conservative scheme \cite{Pimentel2022,Castro2008}
\begin{equation}\label{e:transport_step_scheme_semidiscrete}
    \dot{\Q}_j = - \frac{1}{\Delta x} \left( \matr{A}^-_{j+1} \left( \Q_{j+1} - \Q_{j} \right) + \matr{A}^+_{j} \left( \Q_{j} - \Q_{j-1} \right) \right) =: (\advectQ(\Q))_j,
\end{equation}
where the so-called fluctuations $\matr{A}^{\pm}_{j+1} \left( \Q_{j+1} - \Q_{j} \right)$ can be computed in different ways. 
In this paper, we use the notation from \cite{castro2012} and include the system matrix evaluated at the average of adjacent cells as well as a Lax-Friedrichs type numerical diffusion term for stability as follows
\begin{equation}\label{e:matrix_LF}
    \matr{A}^{\pm}_{j+1} = \frac{1}{2}\left( \matr{A}\left( \frac{\Q_j + \Q_{j+1}}{2} \right) \pm \frac{\Delta x}{\Delta t} \matr{I} \right).
\end{equation}

\begin{remark} 
We note that different numerical schemes are possible that can readily be written in the form of \cref{e:transport_step_scheme}. Examples are the FORCE scheme
\begin{equation}\label{e:matrix_FORCE}
    \matr{A}^{\pm,FORCE}_{j+1} = \frac{1}{2}\left( \matr{A}\left( \frac{\Q_j + \Q_{j+1}}{2} \right) \pm \left(\frac{\Delta x}{2 \Delta t} \matr{I} + \frac{\Delta t}{2 \Delta x} \matr{A}\left( \frac{\Q_j + \Q_{j+1}}{2} \right)^2 \right)\right),
\end{equation}
or the upwind scheme
\begin{equation}\label{e:matrix_upwind}
    \matr{A}^{\pm,upwind}_{j+1} = \frac{1}{2}\left( \matr{A}\left( \frac{\Q_j + \Q_{j+1}}{2} \right) \pm \left|\matr{A}\left( \frac{\Q_j + \Q_{j+1}}{2} \right)\right|\right).
\end{equation}
However, these schemes might include higher computational cost. For details we refer to \cite{koellermeier2021d}.
\end{remark}

Using an explicit Euler time integration scheme for \cref{e:transport_step_scheme_semidiscrete} yields the time discrete update formula of the transport step
\begin{equation}\label{e:transport_step_scheme}
    \Q^{n+1}_j = \Q^{n}_j - \frac{\Delta t}{\Delta x} \left( \matr{A}^-_{j+1} \left( \Q^n_{j+1} - \Q^n_{j} \right) + \matr{A}^+_{j} \left( \Q^n_{j} - \Q^n_{j-1} \right) \right).
\end{equation}

Based on the decomposition of variables and the system matrix given in \cref{eq-def:statematrixV} and \cref{e:block_matrices}, a numerical scheme can be formulated.

For the first order Lax-Friedrichs scheme used in \cref{e:matrix_LF} the fluctuations for both sets of variables can then be derived as follows.
\begin{align}
    \matr{A}^{\pm}_{j+1} =& \frac{1}{2}\left( \matr{A}\left( \frac{\Q_j + \Q_{j+1}}{2} \right) \pm \frac{\Delta x}{\Delta t} \matr{I} \right)\\
    =& \frac{1}{2}\left(
        \begin{bmatrix}
           \matr{A}_{\vect{u}\vect{u}}\left( \frac{\U_j + \U_{j+1}}{2},\frac{\V_j + \V_{j+1}}{2} \right) & \matr{A}_{\vect{u}\vect{v}}\left( \frac{\U_j + \U_{j+1}}{2},\frac{\V_j + \V_{j+1}}{2} \right) \\ \matr{A}_{\vect{v}\vect{u}}\left( \frac{\U_j + \U_{j+1}}{2},\frac{\V_j + \V_{j+1}}{2} \right) & \matr{A}_{\vect{v}\vect{v}}\left( \frac{\U_j + \U_{j+1}}{2},\frac{\V_j + \V_{j+1}}{2} \right)
        \end{bmatrix} \pm \frac{\Delta x}{\Delta t}
        \begin{bmatrix}
           \matr{I}_{2} & \matr{0} \\ \matr{0} & \matr{I}_{N}
        \end{bmatrix} \right)\\
    =& \frac{1}{2}
        \begin{bmatrix}
           \matr{A}_{\vect{u}\vect{u}}\left( \frac{\U_j + \U_{j+1}}{2},\frac{\V_j + \V_{j+1}}{2} \right) \pm \frac{\Delta x}{\Delta t} \matr{I}_{2}  & \matr{A}_{\vect{u}\vect{v}}\left( \frac{\U_j + \U_{j+1}}{2},\frac{\V_j + \V_{j+1}}{2} \right) \\ \matr{A}_{\vect{v}\vect{u}}\left( \frac{\U_j + \U_{j+1}}{2},\frac{\V_j + \V_{j+1}}{2} \right) & \matr{A}_{\vect{v}\vect{v}}\left( \frac{\U_j + \U_{j+1}}{2},\frac{\V_j + \V_{j+1}}{2} \right) \pm \frac{\Delta x}{\Delta t} \matr{I}_{N}
        \end{bmatrix} \\
    =& \frac{1}{2}
        \begin{bmatrix}
           \matr{A}_{j+1,\vect{u}\vect{u}} \pm \frac{\Delta x}{\Delta t} \matr{I}_{2}  & \matr{A}_{j+1,\vect{u}\vect{v}} \\ \matr{A}_{j+1,\vect{v}\vect{u}} & \matr{A}_{j+1,\vect{v}\vect{v}} \pm \frac{\Delta x}{\Delta t} \matr{I}_{N}
        \end{bmatrix}\\
    =& \frac{1}{2}
        \begin{bmatrix}
           \matr{A}^{\pm}_{j+1,\vect{u}\vect{u}} & \matr{A}^{\pm}_{j+1,\vect{u}\vect{v}} \\ \matr{A}^{\pm}_{j+1,\vect{v}\vect{u}} & \matr{A}^{\pm}_{j+1,\vect{v}\vect{v}}
        \end{bmatrix},
\end{align}
where we used the abbreviations $\matr{A}_{j+1,\vect{u}\vect{u}} := \matr{A}_{\vect{u}\vect{u}}\left( \frac{\U_j + \U_{j+1}}{2},\frac{\V_j + \V_{j+1}}{2} \right)$ as well as the identity matrices $\matr{I}_{2} \in \mathbb{R}^{2 \times 2}$ and $\matr{I}_{N} \in \mathbb{R}^{N \times N}$.

The semi-discrete version of \cref{e:transport_step_scheme_semidiscrete} then reads
\begin{align}
    \dot{\Q}_j =& - \frac{1}{\Delta x} \left( \matr{A}^-_{j+1} \left( \Q_{j+1} - \Q_{j} \right) + \matr{A}^+_{j} \left( \Q_{j} - \Q_{j-1} \right) \right)\\
    \begin{bmatrix}
        \dot{\U_j} \\ \dot{\V_j}
    \end{bmatrix}
    =& - \frac{1}{\Delta x} \left( \frac{1}{2}
    \begin{bmatrix}
        \matr{A}^{-}_{j+1,\vect{u}\vect{u}} & \matr{A}^{-}_{j+1,\vect{u}\vect{v}} \\ \matr{A}^{-}_{j+1,\vect{v}\vect{u}} & \matr{A}^{-}_{j+1,\vect{v}\vect{v}}
    \end{bmatrix}
    \begin{bmatrix}
        \U_{j+1} - \U_{j} \\ \V_{j+1} - \V_{j}
    \end{bmatrix} +
    \begin{bmatrix}
        \matr{A}^{+}_{j,\vect{u}\vect{u}} & \matr{A}^{+}_{j,\vect{u}\vect{v}} \\ \matr{A}^{+}_{j,\vect{v}\vect{u}} & \matr{A}^{+}_{j,\vect{v}\vect{v}}
    \end{bmatrix}
    \begin{bmatrix}
        \U_{j} - \U_{j-1} \\ \V_{j} - \V_{j-1}
    \end{bmatrix} \right)\\
    =& \, \frac{1}{2 \Delta t}
    \begin{bmatrix}
        \U_{j+1} - 2 \U_{j} + \U_{j-1} \\ \V_{j+1} - 2 \V_{j} + \V_{j-1}
    \end{bmatrix} + \frac{1}{2 \Delta x} \cdot \dots\\
    &
    \begin{bmatrix}
        \matr{A}_{j+1,\vect{u}\vect{u}} \left(\U_{j+1} - \U_{j}\right) + \matr{A}_{j+1,\vect{u}\vect{v}} \left(\V_{j+1} - \V_{j}\right) + \matr{A}_{j,\vect{u}\vect{u}} \left(\U_{j} - \U_{j-1}\right) + \matr{A}_{j,\vect{u}\vect{v}} \left(\V_{j} - \V_{j-1}\right)  \\ \matr{A}_{j+1,\vect{v}\vect{u}} \left(\U_{j+1} - \U_{j}\right) + \matr{A}_{j+1,\vect{v}\vect{v}} \left(\V_{j+1} - \V_{j}\right) + \matr{A}_{j,\vect{v}\vect{u}} \left(\U_{j} - \U_{j-1}\right) + \matr{A}_{j,\vect{v}\vect{v}} \left(\V_{j} - \V_{j-1}\right)
    \end{bmatrix}\\
     =:& \begin{bmatrix} (\advectU(\U, \V))_j \\ (\advectV(\U, \V))_j \end{bmatrix}.
     \label{eq:RHSV-spacediscrete}
\end{align}
The matrices $\matr{A}_{j+1,\vect{u}\vect{u}}, \matr{A}_{j+1,\vect{u}\vect{v}}, \matr{A}_{j+1,\vect{v}\vect{u}}, \matr{A}_{j+1,\vect{v}\vect{v}}$ are given by \cref{e:block_matrices} evaluated at the averaged variables as follows
\begin{align}
    \matr{A}_{j+1,\vect{u}\vect{u}} =
        \begin{bmatrix}
             & 1 \\
            g h-\left(\frac{1}{2}\left(u_{m,j+1}+u_{m,j}\right)\right)^2-\frac{1}{3}\left(\frac{1}{2}\left(\alpha_{1,j+1}+\alpha_{1,j}\right)\right)^2 & u_{m,j+1}+u_{m,j}
        \end{bmatrix}\in \mathbb{R}^{2 \times 2},
\end{align}
\begin{align}
    \matr{A}_{j+1,\vect{u}\vect{v}} =
        \begin{bmatrix}
             & & \qquad\qquad\qquad\qquad\quad\,\,\,\\
            \frac{1}{3} \left(\alpha_{1,j+1}+\alpha_{1,j}\right) & & \quad
        \end{bmatrix}\in \mathbb{R}^{2 \times N},
\end{align}
\begin{align}
    \matr{A}_{j+1,\vect{v}\vect{u}} =
        \begin{bmatrix}
            -\frac{1}{2}\left(u_{m,j+1}+u_{m,j}\right)\left(\alpha_{1,j+1}+\alpha_{1,j}\right) & \alpha_{1,j+1}+\alpha_{1,j} \\
            -\frac{1}{6} \left(\alpha_{1,j+1}+\alpha_{1,j}\right)^2 & \\
            &                        \\
            &                        \\
            &
        \end{bmatrix} \in \mathbb{R}^{N \times 2},
\end{align}
\begin{align}
    \matr{A}_{j+1,\vect{v}\vect{v}} = \left(\frac12(\alpha_{1,j+1}^{}+\alpha_{1,j}^{})\mathbf{A}+\frac12(u_{m,j+1}^{}+u_{m,j}^{})\mathbf{I}_N\right) \in \mathbb{R}^{N \times N},
\end{align}
where the matrix $\mathbf{A}\in\mathbb{R}^{N\times N}$ is zero except for the off-diagonal entries
\begin{align}
    A_{j,j+1} = \frac{j+2}{2j+3}, \qquad A_{j,j-1} = \frac{j-1}{2j-1}.
\end{align}

The semi-discrete transport step using an explicit Euler time integration scheme for \cref{eq:RHSV-spacediscrete} leads to the time discrete transport update for the first two equations (macro transport step)
\begin{align}
    \U_j^{n+1} =& \frac{1}{2} \left( \U_{j+1}^n + \U_{j-1}^n \right) + \frac{\Delta t}{2 \Delta x} \left( \matr{A}_{j+1,\vect{u}\vect{u}} \left(\U_{j+1}^n - \U_{j}^n\right) + \matr{A}_{j+1,\vect{u}\vect{v}} \left(\V_{j+1}^n - \V_{j}^n\right) \right) \\
    & \qquad +\frac{\Delta t}{2 \Delta x} \left( \matr{A}_{j,\vect{u}\vect{u}} \left(\U_{j}^n - \U_{j-1}^n\right) + \matr{A}_{j,\vect{u}\vect{v}} \left(\V_{j}^n - \V_{j-1}^n\right)\right),
\end{align}
and to the time discrete transport update for the last $N$ moment equations (micro transport step)
\begin{align}
    \V_j^{n+1} =& \frac{1}{2} \left( \V_{j+1}^n + \V_{j-1}^n \right) + \frac{\Delta t}{2 \Delta x} \left( \matr{A}_{j+1,\vect{v}\vect{u}} \left(\U_{j+1}^n - \U_{j}^n\right) + \matr{A}_{j+1,\vect{v}\vect{v}} \left(\V_{j+1}^n - \V_{j}^n\right) \right) \\
    & \qquad +\frac{\Delta t}{2 \Delta x} \left( \matr{A}_{j,\vect{v}\vect{u}} \left(\U_{j}^n - \U_{j-1}^n\right) + \matr{A}_{j,\vect{v}\vect{v}} \left(\V_{j}^n - \V_{j-1}^n\right)\right).
\end{align}

Note that the solution of the micro-step for $\V_j^{n+1}$ is later performed with the known values $\widetilde{\U}_j^{n+1}$ from the first macro step, see \cref{e:transport_micro}. The solution of the macro step should be performed with a conservative scheme, to ensure mass conservation.

\end{document}